\pdfoutput=1	% compile with pdflatex (ensures transparency of figures)
\documentclass{article}
\usepackage[dvipsnames]{xcolor}
\usepackage[hidelinks, hypertexnames=false]{hyperref}
\usepackage{arxiv}  % paper style
\usepackage[utf8]{inputenc} % allow utf-8 input
\usepackage[T1]{fontenc}    % use 8-bit T1 fonts
\usepackage{url}            % simple URL typesetting
\usepackage{booktabs}       % professional-quality tables
\usepackage{nicefrac}       % compact symbols for 1/2, etc.
\usepackage{amsfonts,amssymb,mathtools,amsthm, mathtools}  % for maths
\usepackage{url}
\usepackage{pgfplots}
\pgfplotsset{compat=1.17}
\usepackage{tikz}   % for figures
\usepackage[caption=false]{subfig}
\usepackage{circuitikz} % for arrows in figures
\usepackage{bm, bbm} % bold math font
\usepackage{orcidlink}  % orcid 
\usepackage{dsfont} % for the indicator function symbol
\allowdisplaybreaks % breaking long models on more pages
\usepackage{multirow}   % multirow table
\usepackage{natbib}
\definecolor{mycol1}{HTML}{4F5D75}
\definecolor{mycol2}{HTML}{EF8354}
\definecolor{mycol3}{HTML}{2D3142}
\definecolor{mycol4}{HTML}{650000}
\definecolor{mycol5}{HTML}{BFC0C0}

\newcommand{\name}[1]{\textsf{#1}}

\theoremstyle{definition}
\newtheorem{definition}{Definition}
\newtheorem{example}{Example}

\newcommand\freefootnote[1]{%
	\begin{NoHyper}
		\renewcommand\thefootnote{}\footnote{#1}%
		\addtocounter{footnote}{-1}%
	\end{NoHyper}
}

\title{Multi-Objective Linear Ensembles for Robust and Sparse Training of Few-Bit Neural Networks}

\author{
	Ambrogio Maria Bernardelli$^{*}$ \orcidlink{0000-0002-2328-7062} \\ \texttt{ambrogiomaria.bernardelli01@universitadipavia.it} \And
	Stefano Gualandi$^{*}$ \orcidlink{0000-0002-2111-3528} \\ \texttt{stefano.gualandi@unipv.it} \And
	Simone Milanesi$^{*}$ \orcidlink{0000-0002-6314-1965} \\ \texttt{simone.milanesi01@universitadipavia.it} \And
	Hoong Chuin Lau$^{\dagger}$ \orcidlink{0000-0002-5326-411X} \\ \texttt{hclau@smu.edu.sg} \And
	Neil Yorke-Smith$^{\ddagger}$ \orcidlink{0000-0002-1814-3515} \\ \texttt{n.yorke-smith@tudelft.nl}
}
\date{\today}

\begin{document}
	\maketitle
	\begin{abstract}
		\noindent Training neural networks using combinatorial optimization solvers has gained attention in recent years.  In low-data settings, the use of state-of-the-art mixed integer linear programming solvers, for instance, has the potential to train exactly a neural network (NN), while avoiding intensive GPU-based training and hyper-parameter tuning, and simultaneously training and sparsifying the network.
		We study the case of few-bit discrete-valued neural networks, both Binarized Neural Networks (BNNs), whose values are restricted to $\pm 1$, and Integer Neural Networks (INNs), whose values lie in a range $\{-P, \ldots, P\}$. Few-bit NNs receive increasing recognition due to their lightweight architecture and ability to run on low-power devices, for example being implemented using boolean operations.
		This paper proposes new methods to improve the training of BNNs and INNs.  Our contribution is a multi-objective ensemble approach, based on training a single NN for each possible pair of classes and applying a majority voting scheme to predict the final output.
		Our approach results in the training of robust sparsified networks, whose output is not affected by small perturbations on the input, and whose number of active weights is as small as possible.
		We empirically compare this \emph{BeMi} approach to the current state-of-the-art in solver-based NN training, and to traditional gradient-based training, focusing on BNN learning in few-shot contexts.  We compare the benefits and drawbacks of INNs versus BNNs, bringing new light to the distribution of weights over the $\{-P, \ldots, P\}$ interval.  Finally, we compare multi-objective versus single-objective training of INNs, showing that robustness and network simplicity can be acquired simultaneously, thus obtaining better test performances.
		While the previous state-of-the-art approaches achieve an average accuracy of $ 51.1\%$ on the MNIST dataset, the \emph{BeMi} ensemble approach achieves an average accuracy of $68.4\%$ when trained with 10 images per class and $81.8\%$ when trained with 40 images per class, whilst having up to $75.3\%$ NN links removed.
	\end{abstract}

	\freefootnote{$^{*}$ Department of Mathematics, University of Pavia, Italy}
	\freefootnote{$^{\dagger}$ School of Computing and Information Systems, Singapore Management University, Singapore}
	\freefootnote{$^{\ddagger}$ STAR Lab, Delft University of Technology, Delft, The Netherlands}
	
	\newpage
	\section{Introduction}
	\label{sec:intro}
	
	State-of-the-art deep neural networks (NNs) contain a huge number of neurons organized in many layers, and they require an immense amount of data for training \citep{natDL}. 
	The training process is computationally demanding and is typically performed by stochastic gradient descent algorithms running on large GPU- or even TPU-based clusters.
	Whenever the trained (deep) neural network contains many neurons, also the network deployment is computationally demanding.
	However, in some real-life applications, extensive GPU-based training might be infeasible, training data might be scarce with only a few data points per class, or the hardware using the NN at inference time might have a limited computational power, as for instance, whenever the NN is executed on an industrial embedded system \citep{blott}.

	Binarized Neural Networks (BNNs) were introduced in \citet{Hubara} as a response to the challenge of running NNs on low-power devices.
	BNNs contain only binary weights and binary activation functions, and hence they can be implemented using only efficient bit-wise operations.
	However, the training of BNNs raises interesting challenges for gradient-based approaches due to their combinatorial structure.
	In previous works \citep{Icarte}, it is shown that the training of a BNN can be performed by using combinatorial optimization solvers: a hybrid constraint programming (CP) and mixed integer programming (MIP) approach outperformed the stochastic gradient approach proposed by \citet{Hubara} if restricted to a few-shot-learning context \citep{fewshot}.  \citeauthor{Icarte}'s work has been furthered by a number of authors more recently, as we survey in the next section.
	
	Indeed, combinatorial approaches (principally MIP) for training neural networks, both discrete and continuous, have been employed in the literature, as demonstrated in subsequent
	works, e.g., \citep{patil, kurtz}. 
	MIP optimization has been explored in the machine learning community, as in 
	\citet{DBLP:conf/nips/LiSY20, DBLP:conf/icml/HuangFTDS23, DBLP:conf/icml/YeXWWJ23} for instance.
	Various architectures and activation functions have been utilized in these studies.  Solver-based training has the advantage that, in principle, the optimal NN weights can be found for the training data and that network optimization (e.g., pruning) or adversarial hardening can be performed.
	
	The main challenge in training a NN by an exact MIP-based approach is the limited amount of training data that can be used since, otherwise, the size of the optimization model explodes.
	In the recent work of \citep{Yorke}, the combinatorial training idea was further extended to Integer-valued Neural Networks (INN).  Exploiting the flexibility of MIP solvers, the authors were able to (i) minimize the number of neurons during training and (ii) increase the number of data points used during training by introducing a MIP batch training method.
	
	We remark that training a NN with a MIP-based approach is more challenging than solving a verification problem, as in \citet{fischjo,anderson2020strong}, even if the structure of the nonlinear constraints modelling the activation functions is similar.
	In NNs verification \citep{dilkina}, the weights are given as input, while in MIP-based training, the weights are the \emph{decision variables} that must be computed.

	We note several lines of work aiming at producing compact and simple NNs that maintain acceptable accuracy, e.g., in terms of parameter pruning \citep{serra1,serra2,serracpaior,elaraby2023oamip}, loss function improvement \citep{lfi}, gradient approximation \citep{gdap} and network topology structure \citep{nettopstr}.
	
	In the context of MIP-based training and optimization of NNs, this paper proposes new methods to improve the training of BNNs and INNs.  
	In summary, our contributions are (i) the formulation of a MILP model with a multi-objective target that consists of already existing single-objective steps in a lexicorgraphic order, (ii) the implementation of an ensemble of few-bits NNs in which each of them is specialized in a specific classification task, (iii) the proposal of a voting scheme inspired by One-Versus-One (OVO) strategy, tailored specifically for the constructed ensemble of NNs.
	Our computational results using the MNIST and the Fashion-MNIST dataset show that the \emph{BeMi} ensemble permits to use for training up to 40 data points per class, thanks to the fact that the OVO strategy results in smaller MILPs, reaching an average accuracy of 81.8\% for MNIST and 70.7\% for Fashion-MNIST.
	In addition, thanks to the multi-objective function that minimizes the number of links, i.e. the connections between different neurons, up to 75\% of weights are set to zero for MNIST, and up to 48\% for Fashion-MNIST.
	We also perform additional experiments on the Heart Disease dataset, reaching an average accuracy of 78.5\%. 
	A preliminary report of this work appeared at the LION'23 conference \citep{LION}. This paper better motivates and situates our approach in the state of the art, develops our ensemble approach for INNs (not only BNNs), presents more extensive and improved empirical results, and analyses the distribution of INN weights.
	\paragraph{Outline.} 
	The remainder of this paper is as follows. 
	Section~\ref{sec:rw} situates our work in the literature.
	Section~\ref{sec:notation} introduces the notation and defines the problem of training a single INN with the existing MIP-based methods. 
	Section~\ref{sec:bemi} presents {the \emph{BeMi}} ensemble, the majority voting scheme, and the improved MILP model to train a single INN.
	Section~\ref{sec:results} presents the computational results on the MNIST, Fashion-MNIST and Heart Disease datasets.
	Finally, Section~\ref{sec:conc} concludes the paper with a perspective on future work.
	
	%%%%%%%%%%%%%%%%%%%%%%%%%%%%%%%%%%%%%%%%%%%%%%%%%%%%%%%%%%%%%%%%%%%%%%%%%%%%%%%%%%%%%
	\section{Related Works}
	\label{sec:rw}
	Recently, there has been a growing research interest in studying the impact of machine learning on improving traditional operations research methods (e.g., see \cite{bengio2021machine} and \cite{cappart2023combinatorial}), in designing integrated predictive and modelling frameworks, as in \cite{bergman2022janos}, or in embedding pre-trained machine learning model into an optimization problem (e.g., see \cite{lombardi2018boosting,mistry2021mixed,DBLP:conf/nips/TsayKTM21}).
	In this work, we take a different perspective, and we study how an exact MILP solver can be used for training Machine Learning models, more specifically, to train (Binary or Integer) Neural Networks. In the following paragraphs, we first review two recent applications to MILP solvers in the context of neural networks (i.e., weight pruning and NN verification), and later, we review the few works that tackled the challenge of training an NN using an exact solver.
	
	\paragraph{MIP-based Neural Networks Pruning.}
	Recent works have shown interesting results on {\it pruning} a trained neural network using an optimization approach based on the use of MIP solvers \cite{serra1,serra2,serracpaior,elaraby2023oamip,good2022recall}.
	While pruning a trained neural network, the weights are fixed, and the optimization variables represent the decision to keep or remove an existing weight different from zero.
	
	\paragraph{MIP-based Neural Networks Verification.}
	Another successful application of exact solvers in the context of Neural Networks is for tackling the verification problem, that is, to verify under which conditions the accuracy of a given trained neural network does not deteriorate.
	In other words, in NN verification, the optimization problem consist in finding adversial examples using a minimal distortion of the input data.
	The use of MIP solvers for this application was pioneered in \cite{tjeng2018evaluating}, and later studied in several papers, as for instance, \cite{tjeng2018evaluating,anderson2020strong,botoeva2020efficient,fischjo,tjandraatmadja2020convex}.
	For a broader discussion of the use of polyhedral approach to verification, see Section 4 in \cite{huchette2023deep}.
	In NN verification, the weights of the NN are given as input, and the optimization problem consists of assessing how much the input can change without compromising the output of the network.
	Indeed, also verification is a different optimization problem than the exact training the we discuss in the following sections. 

	\paragraph{Exact training of Neural Networks.}
	% .
	The utilization of 
	MIP approaches for training neural networks has already been explored in the literature, primarily in the context of few-shot learning and NNs with low-bit parameters \citep{Icarte, Yorke}. One of the main advantages of these approaches is the ability to simultaneously train and optimize the network architecture.
	While the modelling that defines the structure of the neural network is rigid and heavily relies on the discrete nature of the parameters, the choice of the objective function provides more flexibility and allows for the optimization of various network characteristics. For instance, it enables minimizing the number of connections in a fixed architecture network, thereby promoting lightweight architectures.
	
	From now on, by INN we indicate a general NN whose weights take value in the set $\{-P, \dots , P\}$, where integer $P \geq 1$.  Notice this choice include the special case of BNNs ($P = 1$).  When referring to an INN with $P > 1$, we will write \emph{non-trivial} INN.
	
	By incorporating the power of both CP and MIP, \citet{Icarte}'s study showcased the effectiveness of leveraging combinatorial optimization methods for training BNNs in a few-shot learning context.  Two types of objective functions were selected to drive the optimization process. The first objective function aimed at promoting a lightweight architecture by minimizing the number of non-zero weights. This objective sought to reduce the overall complexity of the neural network, allowing for efficient computation and resource utilization.
	Note that with this choice, pruning is not necessary: the possibility of setting a weight to zero is equivalent to removing the corresponding weight (i.e., connection).
	In contrast, the second objective function focused on enhancing the robustness of the network, particularly in the face of potential noise in the input data. Robustness here refers to the ability of the network to maintain stable and reliable performance even when the input exhibits variations or disturbances. 
	\citet{Yorke}'s research extends \citet{Icarte}'s single-objective approach to the broader class of Integer Neural Networks (INNs). In addition to leveraging the objectives of architectural lightness and robustness, a novel objective aimed at maximizing the number of correctly classified training data instances is introduced. It is worth noting that this type of objective shares some similarities with the goals pursued by gradient descent methods, albeit with a distinct formulation.
	An interesting observation is that this objective formulation remains feasible in practice, ensuring that a valid solution can be obtained. 
	We remark that both papers propose single-objective models that involve training a single neural network to approximate a multi-classification function.
	
	Ensembles of neural networks are well-known to yield more stable predictions and demonstrate superior generalizability compared to single neural network models \citep{wang}. In this paper, we aim to redefine the concept of structured ensemble by composing our ensemble of several networks, each specialized in a distinct task.
	Given a classification task over $k$ classes, the main idea is to train $\frac{k(k-1)}{2}$ INNs, where every single network learns to discriminate only between a given pair of classes.
	When a new data point (e.g., a new image) must be classified, it is first fed into the $\frac{k(k-1)}{2}$ trained INNs, and later, using a Condorcet-inspired majority voting scheme \citep{condAmSoc}, the most frequent class is predicted as output.  This method is similar to and generalizes the Support Vector Machine - One-Versus-One (SVM-OVO) approach \citep{svmovo}, while it has not yet been applied within the context of neural networks, to the best of our knowledge.
	
	For training every single INN, our approach extends the methods introduced in \citet{Icarte} and \citet{Yorke}, described above.

	\section{Few-bit Neural Networks}
	\label{sec:notation}
	
	In this section, we formally define a Binarized Neural Network (BNN) and an Integer Neural Network (INN) using the notation as in \citet{Icarte} and \citet{Yorke}, while, in the next section, we show how to create a structured ensemble of INNs.
	
	\subsection{Binarized Neural Networks}
	The architecture of a BNN is defined by a set of layers $\mathcal{N} = \{N_0,N_1,\dots,N_L\}$, where $N_l = \{1, \dots, n_l\}$, and $n_l$ is the number of neurons in the $l$-th layer. 
	Let the training set be $\mathcal{X} := \{(\bm{x}^1, y^1), \dots, (\bm{x}^t, y^t)\}$, such that $\bm{x}^i \in \mathbb{R}^{n_0}$ and $y^i \in \{-1, +1\}^{n_L}$ for every $i \in T = \{1, 2, \dots, t\}$. 
	The first layer $N_0$ corresponds to the size of the input data points $\bm{x}^k$.
	Regarding $n_L$, we make the following consideration. For a classification problem with $|\mathcal{I}|$ classes, we set \mbox{$n_L \coloneqq \lceil{\log[2]{|\mathcal{I}|}}\rceil$.} For the case $|\mathcal{I}| = 2$, therefore, $n_L$ will be equal to 1. This is consistent with binary classification problems, as the two classes can be represented as $+1$ and $-1$. When $|\mathcal{I}| = 4$, then $n_L$ will be equal to 2, and the four classes will be represented by $(+1,+1), (+1,-1),(-1,+1)$, and $(-1,-1)$.
	When $|\mathcal{I}|$ is a power of 2, the procedure generalizes in an obvious manner. However, when the number of classes is not a power of 2, we still choose $n_L$ as the nearest integer greater than or equal to the base-2 logarithm of that number, with the caveat that a single network may opt not to classify. For example, if $|\mathcal{I}| = 3$, then $n_L=2$, and $(+1,+1)$ will be associated with the first class, $(+1,-1)$ will be associated with the second class, $(-1,+1)$ will be associated with the third class, while $(-1,-1)$ will be interpreted as `unclassified'.

	The link between neuron $i$ in layer $N_{l-1}$ and neuron $j$ in layer $N_l$ is modelled by weight $w_{ilj} \in \{-1,0,+1\}$. 
	Note that the binarized nature is encoded in the $\pm 1$ weights, while when a weight is set to zero, the corresponding link is removed from the network. 
	Hence, during training, we are also optimizing the architecture of the BNN. 
	
	The activation function is the binary function 
	\begin{equation}\label{eq:activation}
		\rho(x) := 2 \cdot \mathbbm{1} (x \ge 0) - 1,    
	\end{equation}
	that is, a sign function reshaped such that it takes $\pm 1$ values.
	Here the indicator function $\mathbbm{1}(p)$ outputs $+1$ if proposition $p$ is verified, and $0$ otherwise.  This choice for the activation function has been made in line with the literature \citep{Hubara}.
	
	In this paper, we aim to build different MILP models for the simultaneous training and optimization of a network architecture.
	To model the activation function \eqref{eq:activation} of the $j$-th neuron of layer $N_l$ for data point $\bm{x}^k$, we introduce a binary variable $u_{lj}^k \in \{0,1\}$ for the indicator function $\mathbbm{1}(p)$.
	To re-scale the value of $u^k_{lj}$ in $\{-1,+1\}$ and model the activation function value, we introduce the auxiliary variable $z^k_{lj}=(2u^k_{lj} - 1)$.
	For the first input layer, we set $z^k_{0j}=x^k_j$; for the last layer, we account in the loss function whether $z^k_{Lh}$ is different from $y^k_h$.
	The definition of the activation function becomes
	\begin{equation*}
		z^k_{lj} = 
		\rho \left( \sum_{i \in N_{l-1}} z^k_{(l-1)i} w_{ilj} \right) \\
		=
		2 \cdot \mathbbm{1} \left(\sum_{i \in N_{l-1}} z^k_{(l-1)i} w_{ilj} \ge 0\right)  - 1
		= 2 u^k_{lj} - 1.    
	\end{equation*}
	Notice that the activation function at layer $N_l$ gives a nonlinear combination of the output of the neurons in the previous layer $N_{l-1}$ and the weights $w_{ilj}$ between the two layers.
	Section~\ref{sub31} shows how to formulate this activation function in terms of mixed integer linear constraints.
	We remark that the modelling we proposed has already been presented in the literature by \citet{Icarte}.
	
	The choice of a family of parameters $W:=\{w_{ilj}\}_{l\in\{1,\dots,L\},i\in N_{l-1},j\in N_l}$ determines the function
	\begin{equation*}
		f_W : \mathbb{R}^{n_0} \to \{\pm 1\}^{n_L}.
	\end{equation*}
	
	The training of a neural network is the process of computing the family $W$ such that $f_W$ classifies correctly both the given training data, that is, $f_W(\bm{x}^i) = y^i$ for $i=1,\dots,t$, and new unlabelled testing data. 
	
	In the training of a BNN, we follow two machine learning principles for generalization: robustness and simplicity \citep{Icarte}. In      doing so, we target two objectives: 
	(i) the resulting function $f_W$ should generalize from the input data and be {\it robust} to noise in the input data;
	(ii) the resulting network should be {\it simple}, that is, with the smallest number of non-zero weights that permit to achieve the best accuracy.
	
	Regarding the robustness objective, there is argument that deep neural networks have inherent robustness because mini-batch stochastic gradient-based methods implicitly guide toward robust solutions \citep{icarte12,icarte13,icarte25}. 
	However, as shown in \citet{Icarte}, this is false for BNNs in a few-shot learning regime.
	On the contrary, MIP-based training with an appropriate objective function can generalize very well \citep{Icarte,Yorke}, but it does not apply to large training datasets, because the size of the MIP training model is proportional to the size of the training dataset.
	
	One possible way to impose robustness in the context of few-shot learning is to maximize the margins of the neurons, that is, fixing one neuron, we aim at finding an ingoing weights configuration such that for every training input, the entry of the activation function evaluated at that neuron is confidently far away from the discontinuity point. 
	Intuitively, neurons with larger margins require larger changes to their inputs and weights before changing their activation values. 
	This choice is also motivated by recent works showing that margins are good predictors for the generalization of deep convolutional NNs \citep{icarte11}.
	
	Regarding the simplicity objective, a significant parameter is the number of connections \citep{pochipesi}.
	The training algorithm should look for a NN fitting the training data while minimizing the number of non-zero weights. 
	This approach can be interpreted as a simultaneous compression during training, and it has been already explored in recent works \citep{molina,serra3}.
	
	\paragraph{MIP-based BNN training.}
	In \citet{Icarte}, two different MIP models are introduced: the \name{Max-Margin}, which aims to train robust BNNs, and the \name{Min-Weight}, which aims to train simple BNNs. 
	These two models are combined with a CP  model into two hybrid methods \texttt{HW} and \texttt{HA} in order to obtain a feasible solution within a fixed time limit.  \citet{Icarte} employ CP because their MIP models do not scale as the number of training data increases. 
	We remark that in that work, the two objectives, robustness and simplicity, are never optimized simultaneously. 
	
	\paragraph{Gradient-based BNN training.}
	In \citet{Hubara}, a gradient descent-based method is proposed, consisting of a local search that changes the weights to minimize a square hinge loss function. 
	Note that a BNN trained with this approach only learns $\pm 1$ weights. 
	An extension of this method that exploits the same loss function but admits zero-value weights, called {\tt GD$_t$}, is proposed in \citet{Icarte}, to facilitate the comparison with the other approaches. 
	
	\subsection{Integer Neural Networks}
	A more general discrete NN can be obtained when the weights of the network lie in the set $\{-P,-P+1,\dots,-1,0,1,\dots,P-1,P\}$, where $P$ is a positive integer. The resulting network is called INN, and by letting $P=1$, we obtain the BNN presented in the previous subsection. 
	The activation function $\rho$ and the binary variables $u_{lj}^k$ are defined as above.
	The principles leading the training are again simplicity and robustness.
	
	An apparent advantage of using more general integer neural networks lies in the fact that the parameters have increased flexibility while still maintaining their discrete nature. Additionally, by appropriately selecting the parameter $P$, one can determine the number of bits used for each parameter. For instance, $P=1$ corresponds to $1$-bit, $P=3$ corresponds to $2$-bit, $P=8$ corresponds to $3$-bit, and in general, $P=2^{n-1}$ corresponds to $n$-bit.
	
	\paragraph{MIP-based INN training.}
	In \citet{Yorke}, three MIP models are proposed in order to train INNs. 
	The first model, \name{Max-Correct}, is based on the idea of maximizing the number of corrected predicted images;
	the second model, \name{Min-Hinge}, is inspired by the squared hinge loss (compare \citet{Hubara});
	the last model, \name{Sat-Margin}, combines aspects of both the first two models. 
	These three models always produce a feasible solution but use the margins only on the neurons of the last level, obtaining, hence, less robust NNs.
	
	\paragraph{Relations to Quantized Neural Networks.}
	By using an exact MIP solver for training Integer NNs, we are dealing directly with the problem of training a quantized neural network, where all the weights are restricted to take values over a small domain, as discussed above.
	For instance, as reviewed in \cite{gholami2022survey}, there is a growing trend in training NNs using floating point numbers in low precision , that is, using only as few as 8 bits per weight (see for example \cite{banner2018scalable}).
	However, most of the work in the ML literature either focuses on the impact of low-precision arithmetic on the computation of the (stochastic) gradient and in the backpropagation algorithm or focuses on how to {\it quantize} a trained NN by minimizing the deterioration in the accuracy.
	In our work, we take a different perspective on quantization methods, since we do not rely on a gradient-based method to train our INN, but we model and directly solve the problem of training the NN using only a restricted number of integer weights, which is called {\it Integer-only Quantization} in \cite{gholami2022survey}.
	Moreover, by using an exact MIP solver, we can directly find the optimal weights of our (small) INN without running the risks to be trapped into a local minima as stochastic gradient-based methods.

	\section{The \emph{BeMi} ensemble} \label{sec:bemi}
	
	This section first introduces a multi-objective model that allows a simultaneous training and optimization for an INN %(or a BNN) 
	(Section~\ref{sub31}), and then proposes a method for combining a set of neural networks for classification purposes (Section~\ref{sub32}).
	
	\subsection{A multi-objective MILP model for training INNs}
	\label{sub31} 
	For ease of notation, we denote with $\mathcal{L} := \{1, \dots, L\}$ the set of layers and with $\mathcal{L}_2 := \{2, \dots, L\}$, $\mathcal{L}^{L-1} := \{1, \dots, L-1\}$ two of its subsets. 
	We also denote with $\mathfrak{b} \coloneqq \max_{k \in T , j \in N_0  }\{|x^k_j|\} $ a bound on the values of the training data.
	\paragraph{The multi-objective target.}
	A few MIP models are proposed in the literature to train INNs efficiently. 
	In this work, to train a single INN, we use a lexicographic multi-objective function that results in the sequential solution of three different state-of-the-art MIP models: the \name{Sat-Margin} (\texttt{SM}) described in \citet{Yorke}, the \name{Max-Margin} (\texttt{MM}), and the \name{Min-Weight} (\texttt{MW}), both described in \citet{Icarte}. 
	The first model \texttt{SM} maximizes the number of confidently correctly predicted data. 
	The other two models, \texttt{MM} and \texttt{MW}, aim to train a INN following two principles: robustness and simplicity. 
	Our model is based on a lexicographic multi-objective function: first, we train a INN with the model \texttt{SM}, which is fast to solve and always gives a feasible solution.
	Second, we use this solution as a warm start for the \texttt{MM} model, training the INN only with the images that \texttt{SM} correctly classified.
	Third, we fix the margins found with \texttt{MM}, and minimize the number of active weights with \texttt{MW}, finding the simplest INN with the robustness found by \texttt{MM}.
	
	\paragraph{Problem variables.}
	The critical part of our model is the formulation of the nonlinear activation function \eqref{eq:activation}.
	We use an integer variable $w_{ilj} \in \{-P,-P+1, \dots, P\}$ to represent the weight of the connection between neuron $i \in N_{l-1}$ and neuron $j \in N_l$. 
	Variable $u_{lj}^k$ models the result of the indicator function $\mathbbm{1}(p)$  that appears in the activation function $\rho(\cdot)$ for the training instance $\bm{x}^k$. 
	The neuron activation is actually defined as $2 u_{lj}^k - 1$. 
	We introduce auxiliary variables $c^k_{i l j}$ to represent the products $c^k_{i l j} = (2 u_{lj}^k - 1)w_{ilj}$. 
	Note that, while in the first layer, these variables share the same domain of the inputs, from the second layer on, they take values in $\{-P,-P+1, \dots, P\}$.  
	Finally, the auxiliary variables $\hat{y}^k $ represent a predicted label for the input $\bm{x}^k$, and variable $q^k_j$ are used to take into account the data points correctly classified.
	
	The procedure is designed such that the parameter configuration obtained in the first step is used as a warm start for the (\texttt{MM}). Similarly, the solution of the second step is used as a warm start for the solver to solve (\texttt{MW}). In this case, the margins lose their nature as decision variables and become deterministic constants derived from the solution of the previous step.
	
	\paragraph{Sat-Margin (SM) model.}
	We first train our INN using the following \texttt{SM} model.
	\begin{subequations}
		\begin{align}
			\max \quad & \sum_{k \in T} \sum_{j \in N_L}q^k_j \label{eq:fo_sm}\\
			\mbox{s.t.} \quad
			& q_j^k = 1 \implies  \hat{y}_j^k \cdot y_j^k \geq \frac{1}{2}  & \forall j \in N_L, k \in T, \label{sat_impl_1}\\
			& q_j^k = 0 \implies  \hat{y}_j^k \cdot y_j^k \leq \frac{1}{2} - \hat{\epsilon}   & \forall j \in N_L, k \in T, \label{sat_impl_2}\\
			& \hat{y}_j^k = \frac{2}{P \cdot (n_{L-1}+1)}\sum_{i \in N_{L-1}}c^k_{iLj} & \forall j \in N_L, k \in T, \label{eq:sat_pred}\\
			& u^k_{lj} = 1 \implies \sum_{i \in N_{l-1}}c^k_{ilj} \geq0 & \forall l \in \mathcal{L}^{L-1}, j \in N_l, k \in T, \label{sat_impl_3}\\ 
			& u^k_{lj} = 0 \implies  \sum_{i \in N_{l-1}}c^k_{ilj} \leq - \epsilon  & \forall l \in \mathcal{L}^{L-1}, j \in N_l, k \in T,  \label{sat_impl_4}\\
			& c^k_{i 1 j} = x^k_i \cdot w_{i 1 j} & \forall i \in N_0, j \in N_1, k \in T, \label{eq:sat_trivial_mult}\\
			& c^k_{i l j} = (2 u_{(l-1)j}^k- 1)w_{ilj} & \forall l \in \mathcal{L}_2, i \in N_{l-1}, j \in N_l, k \in T, \label{eq:sat_mult}\\
			& q_j^k \in \{0,1\} & \forall j \in N_L, k \in T, \label{eq:sat_q_domain} \\
			& w_{ilj} \in \left\lbrace -P, -P+1, \dots, P \right\rbrace & \forall l \in \mathcal{L}, i \in N_{l-1}, j \in N_l, \label{eq:sat_first_dom}\\
			& u^k_{lj} \in \left\lbrace 0, 1\right\rbrace & \forall l \in \mathcal{L}^{L-1}, j \in N_l, k \in T,  \\
			& c^k_{i1j} \in \left[- P \cdot \mathfrak{b}, P \cdot \mathfrak{b} \right] & \forall i \in N_0, j \in N_1, k\in T,  \\
			& c^k_{ilj} \in \left\lbrace -P, -P+1, \dots, P\right\rbrace & \forall l \in \mathcal{L}_2, i \in N_{l-1}, j \in N_l, k \in T, \label{eq:sat_last_dom}
		\end{align}
	\end{subequations}
	with $\hat{\epsilon} \coloneqq \frac{\epsilon }{2 P \cdot (n_{L-1}+1)}$.
	The objective function \eqref{eq:fo_sm} maximizes the number of data points that are correctly classified.
	Note that $\epsilon$ is a small quantity standardly used to model strict inequalities.
	The implication constraints \eqref{sat_impl_1} and \eqref{sat_impl_2} and constraints \eqref{eq:sat_pred} are used to link the output $\hat{y}^k_j$ with the corresponding variable $q^k_j$ appearing in the objective function. 
	The implication constraints \eqref{sat_impl_3} and \eqref{sat_impl_4} model the result of the indicator function for the $k$-th input data. 
	The constraints \eqref{eq:sat_trivial_mult} and the bilinear constraints \eqref{eq:sat_mult} propagate the results of the activation functions within the neural network.
	We linearize all these constraints with standard big-M techniques \citep{williams2013model}.
	
	The solution of model \eqref{eq:fo_sm}--\eqref{eq:sat_last_dom} gives us the solution vectors $\mathbf{c}_{\texttt{SM}}$, $\mathbf{u}_{\texttt{SM}}$, $\mathbf{w}_{\texttt{SM}}$, $\mathbf{\hat{y}}_{\texttt{SM}}$, $\mathbf{q}_{\texttt{SM}}$. 
	We then define the set
	\begin{equation}
		\hat{T} = \{k \in T \; | \; {q^k_j}_{\texttt{SM}} = 1, \; \forall j \in N_L\},
	\end{equation}
	of confidently correctly predicted images. 
	We use these images as input for the next \name{Max-Margin} \texttt{MM}, and we use the vector of variables $\mathbf{c}_{\texttt{SM}}, \mathbf{u}_{\texttt{SM}}, \mathbf{w}_{\texttt{SM}}$ to warm start the solution of \texttt{MM}.
	
	\paragraph{Max-Margin (MM) model.}
	The second level of our lexicographic multi-objective model maximizes the overall margins of every single neuron activation, with the ultimate goal of training a robust INN.
	Starting from the model {\tt SM}, we introduce the margin variables $m_{lj}$, and we introduce the following \name{Max-Margin} model.
	\begin{subequations}
		\begin{align}
			\max \quad & \sum_{l \in \mathcal{L}} \sum_{j \in N_l}m_{lj}  \\
			\mbox{s.t.} \quad & \mbox{\eqref{eq:sat_trivial_mult}--\eqref{eq:sat_last_dom}} & \forall k \in \hat{T}, \nonumber\\
			& \sum_{i \in N_{L-1}} y^k_j c^k_{iLj} \geq m_{Lj} & \forall j \in N_L, k \in \hat{T}, \label{eq:mm_1}\\
			& u^k_{lj} = 1 \implies \sum_{i \in N_{l-1}}c^k_{ilj} \geq m_{lj} & \forall l \in \mathcal{L}^{L-1}, j \in N_l, k \in \hat{T}, \label{eq:mm_impl_1}\\ 
			& u^k_{lj} = 0 \implies \sum_{i \in N_{l-1}}c^k_{ilj} \leq - m_{lj} & \forall l \in \mathcal{L}^{L-1}, j \in N_l, k \in \hat{T}, \label{eq:mm_impl_2} \\
			& m_{lj} \geq \epsilon & \forall l \in \mathcal{L}, j \in N_l. 
		\end{align}
	\end{subequations}
	Again, we can linearize constraints \eqref{eq:mm_impl_1} and \eqref{eq:mm_impl_2} with standard big-M constraints. 
	This model gives us the solution vectors $\mathbf{c}_{\texttt{MM}}, \mathbf{u}_{\texttt{MM}}, \mathbf{w}_{\texttt{MM}}, \mathbf{m}_{\texttt{MM}}$. We then evaluate $\mathbf{v}_{\texttt{MM}}$ as
	\begin{equation}
		{v_{ilj}}_{\texttt{MM}} =
		\begin{cases}
			0 \quad \textup{when } \; \; {w_{ilj}}_{\texttt{MM}} = 0, \\
			1 \quad \textup{otherwise}, 
		\end{cases}
		\quad \quad  \forall l \in \mathcal{L}, i \in N_{l-1}, j \in N_l.
	\end{equation}
	
	\paragraph{Min-Weight (MW) model.}
	The third level of our multi-objective function minimizes the overall number of non-zero weights, that is, the connections of the trained INN.
	We introduce the new auxiliary binary variable $v_{ilj}$ to model the presence or absence of the link $w_{ilj}$.
	Starting from the solution of model {\tt MM}, we fix $\hat{\mathbf{m}} = 
	\mathbf{m}_{\texttt{MM}}$, and we pass the solution $\mathbf{c}_{\texttt{MM}}, \mathbf{u}_{\texttt{MM}}, \mathbf{w}_{\texttt{MM}}, \mathbf{v}_{\texttt{MM}}$ as a warm start to the following \texttt{MW} model:
	\begin{subequations}
		\begin{align}
			\min \quad & \sum_{l \in \mathcal{L}} \sum_{i \in N_{l-1}}\sum_{j \in N_l}v_{ilj} & \\    
			\mbox{s.t.} \quad    & \mbox{\eqref{eq:sat_trivial_mult}--\eqref{eq:sat_last_dom}} & \forall k \in \hat{T}, \nonumber\\
			& \sum_{i \in N_{L-1}} y^k_j c^k_{iLj} \geq \hat{m}_{Lj} & \forall j \in N_L, k \in \hat{T}, \\
			&  u^k_{lj} = 1 \implies \sum_{i \in N_{l-1}}c^k_{ilj} \geq \hat{m}_{lj} & \forall l \in \mathcal{L}^{L-1}, j \in N_l, k \in \hat{T}, \\ \label{neg_imp_w}
			& u^k_{lj} = 0 \implies \sum_{i \in N_{l-1}}c^k_{ilj} \leq  - \hat{m}_{lj} & \forall l \in \mathcal{L}^{L-1}, j \in N_l, k \in \hat{T},  \\
			&  -v_{ilj} \cdot P \leq w_{ilj} \leq v_{ilj} \cdot P & \forall l \in \mathcal{L}, i \in N_{l-1}, j \in N_l, \label{eq:abs_weight}\\
			& v_{ilj} \in \left\lbrace 0, 1 \right\rbrace & \forall l \in \mathcal{L}, i \in N_{l-1}, j \in N_l. 
		\end{align}
	\end{subequations}
	Note that whenever $v_{ilj}$ is equal to zero, the corresponding weight $w_{ilj}$ is set to zero due to constraint \eqref{eq:abs_weight}, and, hence, the corresponding link can be removed from the network.
	\paragraph{Lexicographic multi-objective.}
	By solving the three models {\tt SM}, {\tt MM}, and {\tt MW}, sequentially, we first maximize the number of input data that is correctly classified, then we maximize the margin of every activation function, and finally, we minimize the number of non-zero weights.
	The solution of the decision variables $w_{ilj}$ of the last model {\tt MW} defines our classification function $f_W: \mathbb{R}^{n_0} \to \{\pm 1\}^{n_L}$.

	\subsection{The \emph{BeMi} structure}
	\label{sub32}
	
	Having explained the various MIP models of INNs, we next introduce our ensemble approach for MIP-based training of INNs.
	
	\paragraph{Ensemble.}
	Define $\mathcal{P} \coloneqq \{\{i, j\} \text{ s.t. } i \neq j, \, i, j \in \mathcal{I}\}$ as the set of all the subsets of the set $\mathcal{I}$ that have cardinality $2$, where $\mathcal{I}$ is the set of the classes of the classification problem.  Then our structured ensemble is constructed in the following way.
	\begin{enumerate}
		\item We train a INN denoted by $\mathcal{N}_{ij}$ for every $\{i,j\} \in \mathcal{P}$, i.e. for each possible pair of elements of $\mathcal{I}$.
		\item When testing a data point, we feed it to our list of trained INNs obtaining a list of predicted labels, namely we obtain the predicted label $\mathfrak{e}_{ij}$ from the network $\mathcal{N}_{ij}$.
		\item We then apply a majority voting system.
	\end{enumerate} 
	
	\smallskip
	The idea behind this structured ensemble is that, given an input $\bm{x}^k$ labelled $l$ $(=y^k)$, the input is fed into $\binom{n}{2}$ networks where $n-1$ of them are trained to recognize an input with label $l$. 
	If all of the networks correctly classify the input $\bm{x}^k$ as $l$, then at most $n-2$ other networks can classify the input with a different label $l'\neq l$, and so the input is correctly labelled with the most occurring label $l$.
	With this approach, if we plan to use $r \in \mathbb{N}$ inputs for each label, we are feeding each of our INNs a total of $2 \cdot r$ inputs instead of feeding $n \cdot r$ inputs to a single large INN. Clearly, when training the networks $\mathcal{N}_{ij}$ and the network $\mathcal{N}_{ik}$, the inputs of the class $i$ are the same, so we only need a total of $r$ inputs for each class. 
	When $n \gg 2$, it is much easier to train our structured ensemble of INNs rather than training one large INN because of the fact that the MILP model size depends linearly on the number of input data.

	\paragraph{Majority voting system.}
	After the training, we feed one input $\bm{x}^k$ to our list of INNs, and we need to elaborate on the set of outputs.
	
	\begin{definition}[Dominant label]
		For every $b \in \mathcal{I}$, we define
		\begin{equation*}
			C_b = \{\{i,j\} \in \mathcal{P} \; | \;  \mathfrak{e}_{ij} = b\},
		\end{equation*}
		and we say that a label $b$ is a \emph{dominant label} if  $|C_b| \geq |C_l|$ for every $l \in \mathcal{I}$. 
		We then define the set of dominant labels
		\begin{equation*}
			\mathcal{D} \coloneqq \{b \in \mathcal{I} \; | \; b \text{ is a dominant label}\}.
		\end{equation*}
	\end{definition}
	Using this definition, we can have three possible outcomes.

	\begin{itemize}
		\item[(a)] There exists a label $i \in \mathcal{I}$ such that $\mathcal{D}= \{i\}$ $\implies$ our input is labelled as $i$.
		\item[(b)] There exist $j, k \in \mathcal{I}$  such that $\mathcal{D} = \{j,k\}$ $\implies$ our input is labelled as $\mathfrak{e}_{jk}$.  
		\item[(c)] There exist more than two dominant labels $\implies$ our input is not classified.
	\end{itemize}
	While case (a) is straightforward, we can label our input even when we do not have a clear winner, that is, when we have trained a INN on the set of labels that are the most frequent (i.e., case (b)).
	Note that the proposed structured ensemble alongside its voting scheme can also be exploited for regular NNs.
	
	\begin{definition}[Label statuses]\label{def:labels}
		In our labelling system, when testing an input, seven different cases (herein called \emph{label statuses}) can arise. The statuses names are of the form `number of the dominant labels + fairness of the prediction'. The first parameter can be $1, 2$, or $o$, where $o$ means `other cases'. The fairness of the prediction is $C$ when it is correct, or $I$ when it is incorrect. The superscripts related to $I'$ and $I''$ only distinguish between different cases. These cases are described through the following tree diagram. \vspace*{15pt}
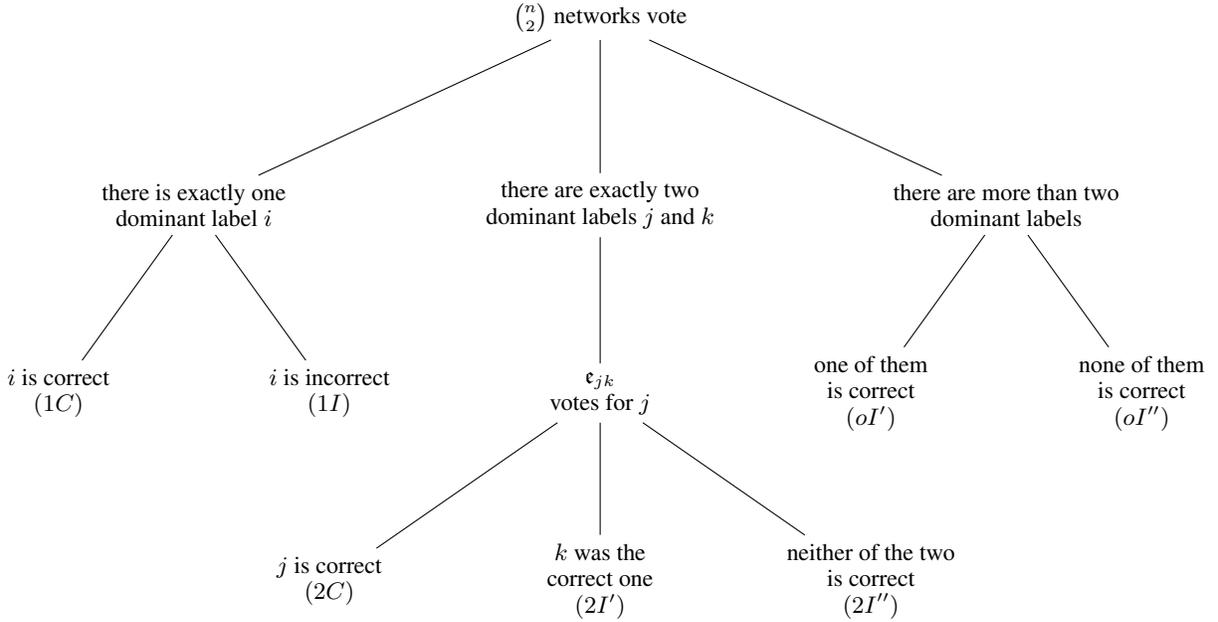
		\begin{center}
			{\small
				\begin{tikzpicture}[black, level distance=2.5cm,
					level 1/.style={sibling distance=5.4cm},
					level 2/.style={sibling distance=3.6cm}]
					\node { $\binom{n}{2}$ networks vote}
					child {node[align=center] {there is exactly one\\ dominant label $i$}
						child {node[align=center] {$i$ is correct \\ $(1C)$ }}
						child {node[align=center] {$i$ is incorrect \\ $(1I)$}}}
					child {node[align=center] {there are exactly two \\ dominant labels $j$ and $k$}
						child {node[align=center] {$ \mathfrak{e}_{jk}$ \\ votes for $j$}
							child {node[align=center] {$j$ is correct \\ $(2C)$}}
							child {node[align=center] {$k$ was the\\  correct one \\ $(2I')$}}
							child {node[align=center] {neither of the two\\ is correct\\ $(2I'')$}}}}
					child {node[align=center] {there are more than two \\ dominant labels}
						child {node[align=center] {one of them \\is correct \\ \textbf{$(oI')$}}}
						child {node[align=center] {none of them \\is correct \\ \textbf{$(oI'')$}}}};
				\end{tikzpicture}
			}
		\end{center}
		
	\end{definition}
	
	The cases in which the classification algorithm classifies correctly are therefore only $(1C)$ and $(2C)$.
	Note that every input test will fall into exactly one label status.

	\smallskip
	\begin{example}
		Let us take $\mathcal{I} = \{bird,\, cat,\, dog,\, frog\}$.
		Note that, in this case, we have to train $\binom{4}{2} = 6$ networks:
		\begin{equation*}
			\mathcal{N}_{\{bird, \, cat\}}, \; \;  \mathcal{N}_{\{bird,\, dog\}}, \; \; \mathcal{N}_{\{bird,\, frog\}}, \; \;
			\mathcal{N}_{\{cat,\, dog\}}, \; \; \mathcal{N}_{\{cat,\, frog\}}, \; \;
			\mathcal{N}_{\{dog,\, frog\}},
		\end{equation*}
		the first one distinguishes between $bird$ and $cat$, the second one between $bird$ and $dog$, and so on.  A first input could have the following predicted labels:
		\begin{align*}
			\mathfrak{e}_{\{bird, \, cat\}} &= bird, &  \mathfrak{e}_{\{bird,\, dog\}} & = bird, & \mathfrak{e}_{\{bird,\, frog\}} &= frog, \\ 
			\mathfrak{e}_{\{cat, \, dog\}}  &= cat,  & 
			\mathfrak{e}_{\{cat, \, frog\}}  &= cat, & 
			\mathfrak{e}_{\{dog, \, frog\}}  &= dog.
		\end{align*}
		We would then have 
		\begin{align*}
			C_{bird} & = \{\{bird, \, cat\}, \{bird, \, dog\}\}, & C_{cat} & = \{\{cat, \, dog\}, \{cat, \, frog\}\}, \\
			C_{dog} & = \{\{dog, \, frog\}\}, & C_{frog} &= \{\{bird, \, frog\}\}.  
		\end{align*}
		In this case $\mathcal{D} = \{bird, \, cat\}$ because $|C_{bird}| = |C_{cat}| = 2 > 1 =   |C_{dog}| = |C_{frog}|$ and we do not have a clear winner, but since $|\mathcal{D}| = 2$, we have trained a network that distinguishes between the two most voted labels, and so we use its output as our final predicted label, labelling our input as $\mathfrak{e}_{\{bird, \, cat\}} = bird$. If $bird$ is the right label we are in label status \textup{$(2C)$}, if the correct label is $cat$, we are in label status \textup{$(2I')$}. Else we are in label status \textup{$(2I'')$}.
	\end{example}
	
	\smallskip
	\begin{example}
		Let us take $\mathcal{I} = \{0,1,\dots, 9\}$. 
		Note that, in this case, we have to train $\binom{10}{2} = 45$ networks and that $|C_b| \leq 9$ for all $b \in \mathcal{I}$. 
		Hence, an input could be labelled as follows:
		\begin{align*}
			& C_0 = (\{0, i\})_{i = 1,2,3,5,7,8}\; , & &  C_1 = (\{1, i\})_{i = 5, 6} \; ,& & C_2 = (\{2, i\})_{i=1,5,8} \; ,\\
			&C_3 =(\{3, i\})_{i = 1,2,4,5} \; , && C_4 = (\{4, i\})_{i=0,1,2,5,6,7,9} \; ,&& C_5 = (\{5, i\})_{i=6,7} \; ,\\
			& C_6 = (\{6, i\})_{i=0,2,3,7} \; , && C_7 = (\{7, i\})_{i=1,2,3} \; , && \\ & C_8 = (\{8, i\})_{i=1,3,4,5,6,7} \; ,
			&& C_9 = (\{9, i\})_{i=0,1,2,3,5,6,7,8} \; . &&
		\end{align*}
		Visually, we can represent an input being labelled as above with the following scheme:
		
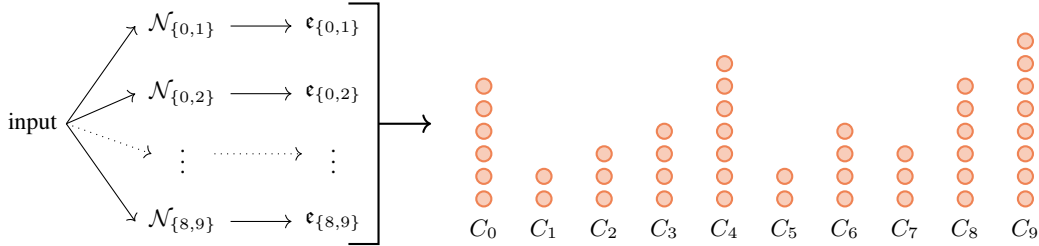
\begin{center}
			\begin{tikzpicture}
				\node[rectangle] (1) at (0, 0){\small \textcolor{black}{input}};
				\node[circle] (2) at (2, 1.3){\small $\mathcal{N}_{\{0, 1\}}$};
				\node[circle] (3) at (2, 0.4){\small $\mathcal{N}_{\{0, 2\}}$};
				\node[circle] (4) at (2, -0.4){\small $\vdots$};
				\node[circle] (5) at (2, -1.3){\small $\mathcal{N}_{\{8, 9\}}$};
				\draw[->] (1.east) -- (2.west);
				\draw[->] (1.east) -- (3.west);
				\draw[->, dotted] (1.east) -- (4.west);
				\draw[->] (1.east) -- (5.west);
				\node[circle] (9) at (4, 1.3){\small $\mathfrak{e}_{\{0, 1\}}$};
				\node[circle] (10) at (4, 0.4){\small $\mathfrak{e}_{\{0, 2\}}$};
				\node[circle] (11) at (4, -0.4){\small $\vdots$};
				\node[circle] (12) at (4, -1.3){\small $\mathfrak{e}_{\{8, 9\}}$};
				\draw[->] (2.east) -- (9.west);
				\draw[->] (3.east) -- (10.west);
				\draw[->, dotted] (4.east) -- (11.west);
				\draw[->] (5.east) -- (12.west);
				\draw[thick] (4.2, 1.6) -- (4.6, 1.6) -- (4.6, -1.6) -- (4.2, -1.6);
				\draw[->, thick] (4.6, 0) -- (5.3, 0);
				\node[circle] (20) at (6, -1.4){\small $C_0$};
				\draw[fill=mycol2, thick, fill opacity=0.4, color=mycol2] (6,-1.0) circle (0.1);
				\draw[fill=mycol2, thick, fill opacity=0.4, color=mycol2] (6,-0.7) circle (0.1);
				\draw[fill=mycol2, thick, fill opacity=0.4, color=mycol2] (6,-0.4) circle (0.1);
				\draw[fill=mycol2, thick, fill opacity=0.4, color=mycol2] (6,-0.1) circle (0.1);
				\draw[fill=mycol2, thick, fill opacity=0.4, color=mycol2] (6,0.2) circle (0.1);
				\draw[fill=mycol2, thick, fill opacity=0.4, color=mycol2] (6,0.5) circle (0.1);
				\node[circle] (19) at (6.8, -1.4){\small $C_1$};
				\draw[fill=mycol2, thick, fill opacity=0.4, color=mycol2] (6.8,-1.0) circle (0.1);
				\draw[fill=mycol2, thick, fill opacity=0.4, color=mycol2] (6.8,-0.7) circle (0.1);
				\node[circle] (21) at (7.6, -1.4){\small $C_2$};
				\draw[fill=mycol2, thick, fill opacity=0.4, color=mycol2] (7.6,-1.0) circle (0.1);
				\draw[fill=mycol2, thick, fill opacity=0.4, color=mycol2] (7.6,-0.7) circle (0.1);
				\draw[fill=mycol2, thick, fill opacity=0.4, color=mycol2] (7.6,-0.4) circle (0.1);
				\node[circle] (22) at (8.4, -1.4){\small $C_3$};
				\draw[fill=mycol2, thick, fill opacity=0.4, color=mycol2] (8.4,-1.0) circle (0.1);
				\draw[fill=mycol2, thick, fill opacity=0.4, color=mycol2] (8.4,-0.7) circle (0.1);
				\draw[fill=mycol2, thick, fill opacity=0.4, color=mycol2] (8.4,-0.4) circle (0.1);
				\draw[fill=mycol2, thick, fill opacity=0.4, color=mycol2] (8.4,-0.1) circle (0.1);
				\node[circle] (21) at (9.2, -1.4){\small $C_4$};
				\draw[fill=mycol2, thick, fill opacity=0.4, color=mycol2] (9.2,-1.0) circle (0.1);
				\draw[fill=mycol2, thick, fill opacity=0.4, color=mycol2] (9.2,-0.7) circle (0.1);
				\draw[fill=mycol2, thick, fill opacity=0.4, color=mycol2] (9.2,-0.4) circle (0.1);
				\draw[fill=mycol2, thick, fill opacity=0.4, color=mycol2] (9.2,-0.1) circle (0.1);
				\draw[fill=mycol2, thick, fill opacity=0.4, color=mycol2] (9.2,0.2) circle (0.1);
				\draw[fill=mycol2, thick, fill opacity=0.4, color=mycol2] (9.2,0.5) circle (0.1);
				\draw[fill=mycol2, thick, fill opacity=0.4, color=mycol2] (9.2,0.8) circle (0.1);
				\node[circle] (21) at (10, -1.4){\small $C_5$};
				\draw[fill=mycol2, thick, fill opacity=0.4, color=mycol2] (10,-1.0) circle (0.1);
				\draw[fill=mycol2, thick, fill opacity=0.4, color=mycol2] (10,-0.7) circle (0.1);
				\node[circle] (21) at (10.8, -1.4){\small $C_6$};
				\draw[fill=mycol2, thick, fill opacity=0.4, color=mycol2] (10.8,-1.0) circle (0.1);
				\draw[fill=mycol2, thick, fill opacity=0.4, color=mycol2] (10.8,-0.7) circle (0.1);
				\draw[fill=mycol2, thick, fill opacity=0.4, color=mycol2] (10.8,-0.4) circle (0.1);
				\draw[fill=mycol2, thick, fill opacity=0.4, color=mycol2] (10.8,-0.1) circle (0.1);
				\node[circle] (21) at (11.6, -1.4){\small $C_7$};
				\draw[fill=mycol2, thick, fill opacity=0.4, color=mycol2] (11.6,-1.0) circle (0.1);
				\draw[fill=mycol2, thick, fill opacity=0.4, color=mycol2] (11.6,-0.7) circle (0.1);
				\draw[fill=mycol2, thick, fill opacity=0.4, color=mycol2] (11.6,-0.4) circle (0.1);
				\node[circle] (21) at (12.4, -1.4){\small $C_8$};
				\draw[fill=mycol2, thick, fill opacity=0.4, color=mycol2] (12.4,-1.0) circle (0.1);
				\draw[fill=mycol2, thick, fill opacity=0.4, color=mycol2] (12.4,-0.7) circle (0.1);
				\draw[fill=mycol2, thick, fill opacity=0.4, color=mycol2] (12.4,-0.4) circle (0.1);
				\draw[fill=mycol2, thick, fill opacity=0.4, color=mycol2] (12.4,-0.1) circle (0.1);
				\draw[fill=mycol2, thick, fill opacity=0.4, color=mycol2] (12.4,0.2) circle (0.1);
				\draw[fill=mycol2, thick, fill opacity=0.4, color=mycol2] (12.4,0.5) circle (0.1);
				\node[circle] (21) at (13.2, -1.4){\small $C_9$};
				\draw[fill=mycol2, thick, fill opacity=0.4, color=mycol2] (13.2,-1.0) circle (0.1);
				\draw[fill=mycol2, thick, fill opacity=0.4, color=mycol2] (13.2,-0.7) circle (0.1);
				\draw[fill=mycol2, thick, fill opacity=0.4, color=mycol2] (13.2,-0.4) circle (0.1);
				\draw[fill=mycol2, thick, fill opacity=0.4, color=mycol2] (13.2,-0.1) circle (0.1);
				\draw[fill=mycol2, thick, fill opacity=0.4, color=mycol2] (13.2,0.2) circle (0.1);
				\draw[fill=mycol2, thick, fill opacity=0.4, color=mycol2] (13.2,0.5) circle (0.1);
				\draw[fill=mycol2, thick, fill opacity=0.4, color=mycol2] (13.2,0.8) circle (0.1);
				\draw[fill=mycol2, thick, fill opacity=0.4, color=mycol2] (13.2,1.1) circle (0.1);
			\end{tikzpicture}
		\end{center}
		where we have omitted the name of each element of the set $C_i$ for simplicity: for example, the dots above $C_1$ represent the sets $\{1, 5\}, \{1, 6\}$. Since $\mathcal{D} = \{9\}$, our input is labelled as $9$.
		If $9$ is the right label, we are in label status $(1C)$, if it is the wrong one, we are in label status $(1I)$. 
		If instead $\hat{C}_j = C_j$, $j=0, \dots, 7$, and
		\begin{align*}
			& \hat{C}_8 = (\{8, i\})_{i=1,3,4,5,6,7, 9} \; , && \hat{C}_9 = (\{9, i\})_{i=0,1,2,3,5,6,7} \; ,
		\end{align*}
		then $|\mathcal{\hat{D}}| = |\{4, 8, 9\}| = 3$, so that our input was labelled as $-1$. 
		If the correct label is $4, 8$ or $9$, we are in label status $(mI')$, else we are in label status $(mI'')$. Lastly, if $\bar{C}_j = C_j$, $j \in \{0,1,2,4,5,6,7,8\}$, and
		\begin{align*}
			& \bar{C}_3 = (\{3, i\})_{i=1,2,4,5,9} \; , && \bar{C}_9  = (\{9, i\})_{i=0,1,2,5,6,7,8}
		\end{align*}
		then $|\mathcal{\bar{D}}| = |\{4, 9\}| = 2 $ and since $\{4, 9\} \in \bar{C}_4$ our input is labelled as $4$. If $4$ is the correct label, we are in label status $(2C)$, if $9$ is the correct label, we are in label status $(2I')$, else we are in label status $(2I'')$.
		Note that, for brevity, in this example we used the notation $
		(\{j, i\})_{i = i_1, \dots, i_n} = \{\{j, i_1\}, \dots, \{j, i_n\}\}$, $ j, i_1, \dots, i_n \in \{0, \dots, 9\}
		$.
		
	\end{example}

	%%%%%%%%%%%%%%%%%%%%%%%%%%%%%%%%%%%%%%%%%%%%%%%%%%%%%%%%%

	\section{Empirical Study}
	\label{sec:results}
	
	Having introduced the \emph{BeMi} approach,
	we now undertake a series of six experiments in order to explore the following questions:
	\begin{enumerate}
		\item[-]{\bf Experiment 1:} What is the impact of a three-fold multi-objective model compared to a two-fold or single objective model?  (Recall Section~\ref{sec:bemi}.)
		
		\item[-]{\bf Experiment 2:} How does the \emph{BeMi} ensemble compare with the previous state-of-the-art MIP models for training BNNs in the context of few-shot learning?
		
		\item[-]{\bf Experiment 3:} How does the \emph{BeMi} ensemble scale with the number of training images, considering two different types of BNNs?
		
		\item[-]{\bf Experiment 4:} How does the \emph{BeMi} ensemble perform on various datasets, comparing the running time, the average gap to the optimal training MILP model, and the percentage of links removed?
		
		\item[-]{\bf Experiment~5:} What are the performance differences between a non-trivial INN and a BNN? Do INN exhibit particular weights distribution characteristics? A state-of-the-art comparison is also provided.  
		
		\item[-]{\bf Experiment~6:} How does the \emph{BeMi} ensemble perform on a continuous, low-dimension dataset, comparing BNNs and non-trivial INNs? Do INN exhibit the same weights distribution characteristics found in Experiment 5?
		
	\end{enumerate}
	\paragraph{Datasets.} 
	Three datasets are adopted for the experiments.  We use first the standard MNIST dataset \citep{mnist} for a fair comparison with the literature, and second the larger Fashion-MNIST dataset \citep{fmnist}.
	For these two MNIST datasets, we test our results on $800$ images for each class in order to have the same amount of test data for every class. Note that the MNIST dataset has $10\,000$ test data but they are not uniformly distributed over the $10$ classes.  
	For each experiment, we report the average over five different samples of images, i.e., we perform five different trainings and we report the average over them, while testing the same images.
	The images are sampled uniformly at random in order to avoid overlapping between different experiments.
	Beyond MNIST, we use the Heart Disease dataset \citep{heart} from the UCI repository.
	Table~\ref{tab:dataset} summarizes the datasets.
	
	\begin{table}[t!]
		\centering
		\small
		\begin{tabular}{c@{\hskip 0.3in}c@{\hskip 0.3in}c@{\hskip 0.3in}c@{\hskip 0.3in}c@{\hskip 0.3in}c}
			\toprule    
			Dataset & Number of classes & Input Dimension & Data Values & \# Training Set  & \# Test Set  \\
			\midrule
			MNIST & $10$ & $28 \times 28$ & Integers & $60\,000$  & $10\,000$  \\
			Fashion-MNIST & $10$ & $28 \times 28$ & Integers & $60\,000$  & $10\,000$   \\
			Heart Disease & $2$ & $13$ & Continuous & $920 - x$  &  $x$ \\
			\bottomrule
		\end{tabular}
		\caption{Details of the different datasets exploited in the experiments.}
		\label{tab:dataset}
	\end{table}
	
	\paragraph{Implementation details.}
	As the solver we use Gurobi version 10.0.1 \citep{gurobi} to solve our MILP models. 
	The solver parameters are left to the default values if not specified otherwise.
	Apart from the first experiment, where we chose to consider every model equally, the fraction of time given to each step of the multi-objective model has been chosen accordingly to the importance of finding a feasible and robust solution.  
	All the MILP experiments were run on an HPC cluster running CentOS, using a single node per experiment.
	Each node has an Intel CPU with 8 physical cores working at 2.1 GHz, and 16 GB of RAM.
	In all of our experiments concerning integer-value datasets, we fix the value $\epsilon = 0.1$. Notice that, because of the integer nature of the weights, of the image of the activation function, and of the data, setting $\epsilon$ equal to any number smaller than $1$ is equivalent. When using continuous-value datasets, we fix the value $\epsilon = 1 \cdot 10^{-6}$ in accordance with the default variable precision tolerance of the Gurobi MILP solver we will use. 
	The source code is available on GitHub \citep{codice-nostro}.

	\paragraph{Time limit management.} Concerning the time limits for the different optimization models, the following choices have been made.
	In Experiment 1 and 6, the time limit is equally distributed between the three models
	to have a fair comparison. In Experiment 2, 3, 4, and 5, the majority of the
	imposed time limit was given to the first two models. The first model ensures
	feasibility of the whole pipeline and maximises the number of correctly classified images in the training phase, and it was considered important
	in the context of few-shot learning, since we do not have lots of images as training inputs. The second model was given a bigger time
	limit too because preliminary results, also shown in previous works, highlight
	the fact that the maximization of the margin ensures a better test accuracy
	with respect to the minimization of the links.
	In addition to this, the overall time limits have been chosen based on two criteria. Where comparisons with the literature are made, the selection ensures a fair comparison. In the remaining cases, the choice of time limit has been empirical, aiming to highlight the algorithm's quality.
	
	\subsection{Experiment 1}
	The goal of the first experiment is to study the impact of the multi-objective model composed by \texttt{SM}, \texttt{MM} and \texttt{MW} with respect to the models composed by \texttt{SM} and \texttt{MM}, the one composed by by \texttt{SM} and \texttt{MW}, and only \texttt{SM}, respectively.
	
	The results refer to a BNN specialised in distinguishing between digits 4 and 9 of the MNIST dataset.  These two digits were chosen because their written form can be quite similar. Indeed, among all ten digits, 4 and 9 are very often mistaken for each other, as it is shown in the confusion matrix in Appendix \ref{app:exp}.
	
	The NN architecture consists of two hidden layers and has $[784,4,4,1]$ neurons. The architecture is chosen to mimic the one used in \citet{Icarte}, that is $[784,16,16,10]$, but with fewer neurons. The number of training images varies between $2$, $6$ and $10$, while the test images are $800$ per digit, so $1600$ in total. The imposed time limit is $30$ minutes, equally distributed in the steps of each model: $30$ minutes for the \texttt{SM}, $15+15$ for the \texttt{SM}+\texttt{MW}, $15+15$ for the \texttt{SM}+\texttt{MM}, and $10+10+10$ for \texttt{SM}+\texttt{MM}+\texttt{MW}.
	
\begin{figure}[t!]
		\centering
		\subfloat[t!][Test accuracy of different models.\label{fig:acc_abl}]{%
			\begin{tikzpicture}
				
				\begin{axis}[
					width=6.5cm,
					height=4cm,
					scale only axis,
					xmin=1, xmax=11,
					xtick={2,6,10},
					xticklabels={\small 2,\small 6,\small 10},
					xmajorgrids,
					ymin=40, ymax=100,
					ytick={40, 50, 60, 70, 80, 90, 100},
					yticklabels={\small 40\%, \small 50\%, \small 60\%, \small 70\%, \small 80\%, \small 90\%, \small 100\%},
					ylabel={\small Test Accuracy},
					xlabel={\small Training Images per Digit},
					legend cell align={left},
					ymajorgrids,
					axis lines*=left,
					line width=1.0pt,
					mark size=3.0pt,
					legend style ={ at={(0,1)},
						anchor=north west, draw=black, 
						fill=white,align=left},
					cycle list name=black white,
					smooth
					]
					
					\addplot[color=mycol3, mark=triangle, dash dot, mark options={solid}] coordinates{
						(2, 52.85)
						(6, 55.25)
						(10, 58.95) 
					};
					\addlegendentry{\small \texttt{SM}};

					\addplot[color=mycol4, mark=star, densely dotted, mark options={solid}] coordinates{
						(2, 58.38)
						(6, 61.70)
						(10, 68.89)
					};
					\addlegendentry{\small \texttt{SM}+\texttt{MW}};
					
					\addplot[color=mycol2, mark=square, dashed, mark options={solid}] coordinates{
						(2, 65.05)
						(6, 66.33)
						(10, 75.12) 
					};
					\addlegendentry{\small \texttt{SM}+\texttt{MM}};
					
					\addplot[color=mycol1, mark=o,dotted, mark options={solid}] coordinates{
						(2, 65.13)
						(6, 66.45)
						(10, 74.70) 
					};
					\addlegendentry{\small \texttt{SM}+\texttt{MM}+\texttt{MW}};
					\fill[mycol3, opacity = 0.1] (2, 50) -- (6, 50) -- (10, 54.13) -- (10, 63.25) -- (6, 60.56) -- (2, 60.75);
					\fill[mycol4, opacity = 0.1] (2, 56.93) -- (6, 56.94) -- (10, 67.56) -- (10, 78.75) -- (6, 64.38) -- (2, 64.68);
					\fill[mycol2, opacity = 0.1] (2, 58.38) -- (6, 58.13) -- (10, 72.38) -- (10, 78.69) -- (6, 71.81) -- (2, 69.50);
					\fill[mycol1, opacity = 0.1] (2, 58.13) -- (6, 57.44) -- (10, 71.63) -- (10, 77.88) -- (6, 72.00) -- (2, 69.75);
					
				\end{axis}
			\end{tikzpicture}%
		}\hfill
		\subfloat[t!][Simplicity of different models.\label{fig:weights_abl}]{%
			\begin{tikzpicture}
				\begin{axis}[
					width=6.5cm,
					height=4cm,
					scale only axis,
					xmin=1, xmax=11,
					xtick={2,6,10},
					xticklabels={\small 2,\small 6,\small 10},
					xmajorgrids,
					ymin=0, ymax=60,
					ytick={0, 10, 20, 30, 40, 50, 60},
					yticklabels={\small 0\%, \small 10\%, \small 20\%, \small 30\%, \small 40\%, \small 50\%, \small 60\%,},
					ylabel={\small Non-zero Weights},
					xlabel={\small Training Images per Digit},
					ymajorgrids,
					axis lines*=left,
					line width=1.0pt,
					mark size=3.0pt,
					legend style ={ at={(1.08,1)}, 
						anchor=north west, %draw=black, 
						fill=white,align=left},
					cycle list name=black white,
					smooth
					]
					\addplot[color=mycol3, mark=triangle, dash dot, mark options={solid}] coordinates{
						(2, 35.35)
						(6, 42.96)
						(10, 49.22) 
					};
					
					\addplot[color=mycol4, mark=star, densely dotted, mark options={solid}] coordinates{
						(2, 0.1)
						(6, 0.1)
						(10, 0.13)
					};
					\addplot[color=mycol2, mark=square, dashed, mark options={solid}] coordinates{
						(2, 35.10)
						(6, 42.63)
						(10, 48.90) 
					};
					
					\addplot[color=mycol1, mark=o,dotted, mark options={solid}] coordinates{
						(2, 23.53)
						(6, 26.88)
						(10, 27.24) 
					};
					
					\fill[mycol3, opacity = 0.1] (2, 31.27) -- (6, 38.69) -- (10, 45.06) -- (10, 54.25) -- (6, 48.67) -- (2, 41.41);
					\fill[mycol4, opacity = 0.1] (2, 0.1) -- (6, 0.1) -- (10, 0.1) -- (10, 0.13) -- (6, 0.13) -- (2, 0.1);
					\fill[mycol2, opacity = 0.1] (2, 31.18) -- (6, 38.31) -- (10, 44.74) -- (10, 54.12) -- (6, 48.42) -- (2, 41.19);
					\fill[mycol1, opacity = 0.1] (2, 18.63) -- (6, 21.01) -- (10, 22.72) -- (10, 33.62) -- (6, 39.59) -- (2, 29.34);
					
				\end{axis}
			\end{tikzpicture}%
		}
		\caption{The \texttt{SM}+\texttt{MM}+\texttt{MW} achieve the same accuracy of the \texttt{SM}+\texttt{MM} model, outperforming the \texttt{SM} model, while having the smallest percentage of non-zero weights, a part from the \texttt{SM}+\texttt{MW} model, which has an almost-zero percentage of non-zero weights, but also a lower accuracy of the models that maximize the margins. The dotted lines represent the average accuracy obtained over 5 instances, while the shaded areas highlight the range between the minimum and maximum values.}
		\label{fig:ablation}
	\end{figure}
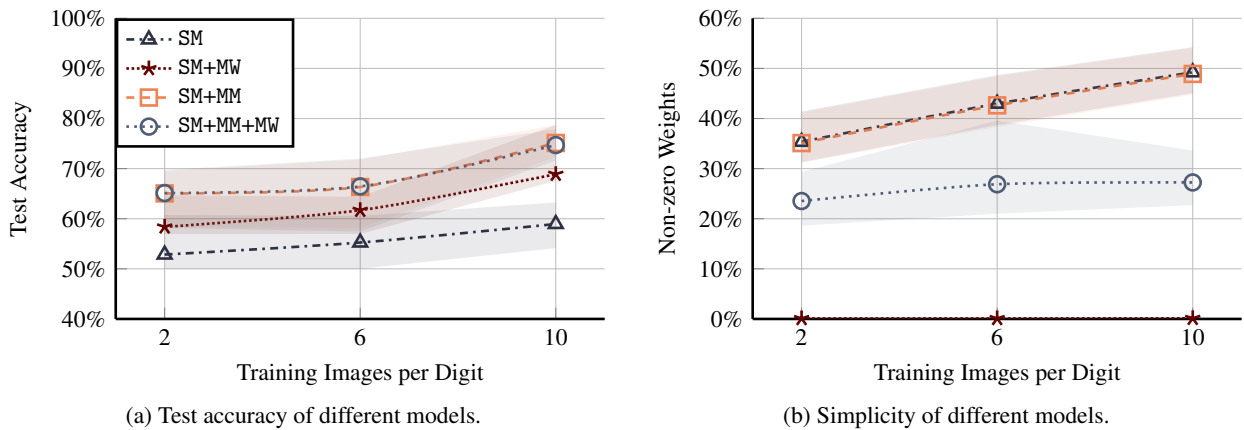

	Figure~\ref{fig:acc_abl} compares the test accuracy of the four hierarchical models, showing how the \texttt{MM} model ensures an increase in accuracy, while the \texttt{MW} allows the network to be pruned without performance being affected. Figure~\ref{fig:weights_abl} displays the percentages of non-zero weights of the three trained models. Note that in this case the training accuracy is always $100\%$ and so we did not add it to the plot. Also, dotted lines represent the average accuracy obtained over 5 instances, while the shaded areas highlight the range between the minimum and maximum values. This will be the case for every other plot if not specified otherwise. The \texttt{MW} step allows the number of non-zero weights to drop without accuracy being affected, hence resulting in an effective pruning. This behaviour is also observed for other couples of digits, even the ones that are easier to distinguish, namely, $1$ and $8$. Results of this experiment are reported in Appendix \ref{app:exp}.
	
	Based on these reasons, for the remaining experiments, we will exclusively employ the model that incorporates all three steps.
	
	\subsection{Experiment 2}
	\begin{figure}[t!]
		\centering
		%%%%%% NEW PICTURE
		\begin{tikzpicture}
			\begin{axis}[
				width=14cm,
				height=6cm,
				scale only axis,
				xmin=0, xmax=11,
				xtick={2,6,10},
				xticklabels={\small 2,\small 6,\small 10},
				xmajorgrids,
				ymin=0, ymax=80,
				ytick={0, 20, 40, 60, 80},
				yticklabels={\small 0\%, \small 20\%, \small 40\%, \small 60\%, \small 80\%},
				ylabel={\small Test Accuracy},
				xlabel={\small Training Images per Digit},
				legend cell align={left},
				ymajorgrids,
				axis lines*=left,
				line width=1.0pt,
				mark size=3.0pt,
				legend style ={ at={(0,1)}, 
					anchor=north west, %draw=black, 
					fill=white,align=left},
				cycle list name=black white,
				smooth
				]
				
				\fill[mycol1, opacity = 0.1] (2, 45.64) -- (6, 60.10) -- (10, 65.49) -- (10, 74.54) -- (6, 67.41) -- (2, 51.40);

				\addplot[color=mycol1, mark=o, dotted, mark options={solid}] coordinates{
					(2, 49.52)
					(6, 63.74)
					(10, 68.35) 
				};
				\addlegendentry{\small \emph{BeMi}};
				\addplot[color=mycol2, mark=diamond, mark options={solid}] coordinates{
					(10, 53.90) 
				};
				\addlegendentry{\small \texttt{MH}};

				\addplot[color=mycol4, mark=star, dashed, mark options={solid}] coordinates{
					(2, 26.53)
					(6, 41.53)
					(10, 52.64) 
				};
				\addlegendentry{\small \texttt{HA}};
				
				\addplot[color=mycol3, mark=square, mark options={solid}] coordinates{
					(1, 12.22) 
				};
				\addlegendentry{\small \texttt{MIP}};
				
				\addplot[color=mycol5, mark=triangle, dash dot, mark options={solid}] coordinates{
					(2, 7.26)
					(6, 13.47)
					(10, 13.61) 
				};
				\addlegendentry{\small \texttt{GD}$_t$};
				
			\end{axis}
		\end{tikzpicture}%

		\caption{Comparison of published approaches vs \emph{BeMi}, in terms of accuracy over the MNIST dataset using few-shot learning with 2, 6, and 10 images per digit.}    
		\label{fig:comparison}
	\end{figure}
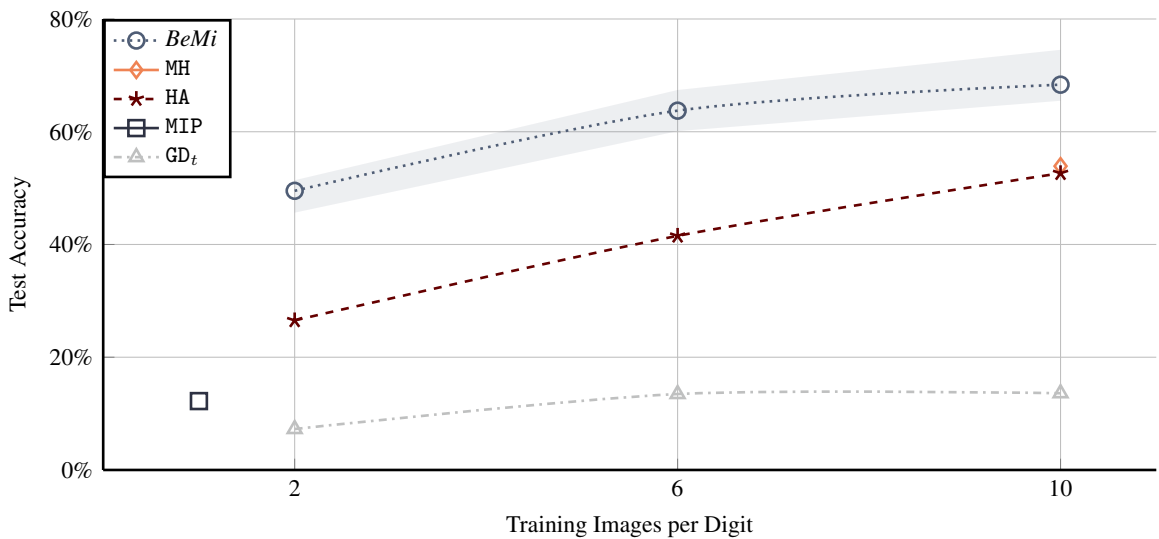
	
	The goal of the second experiment is to compare the \emph{BeMi} ensemble with the following state-of-the-art methods:
	the pure MIP model in \citet{Icarte};
	the hybrid CP--MIP model based on \name{Max-Margin} optimization (\texttt{HA}) \citet{Icarte};
	the gradient-based method {\tt GD}$_t$ introduced in \citet{Hubara} and adapted in \citet{Icarte} to deal with link removal;
	and the \name{Min-hinge} (\texttt{MH}) model proposed in \citet{Yorke}. 
	For the comparison, we fix the setting of \citet{Icarte}, which takes from the MNIST up to 10 images for each class, for a total of 100 training data points, and which uses a time limit of $7\,200$ seconds to solve their MIP training models.
	
	In our experiments, we train the \emph{BeMi} ensemble with $2$, $6$ and $10$ samples for each digit.
	Since our ensemble has 45 BNNs, we leave for the training of each single BNN a maximum of 160 seconds (since $160 \cdot 45 = 7\,200$).
	In particular, we give a $75$ seconds time limit to the solution of \texttt{SM}, $75$ seconds to \texttt{MM}, and $10$ seconds to \texttt{MW}. 
	In all of our experiments, whenever the optimum is reached within the time limit, the remaining time is added to the time limit of the subsequent model. 
	We remark that our networks could be trained in parallel, which would highly reduce the wall-clock runtime. 
	For the sake the completeness, we note that we are using $45 \cdot (784 \cdot 4 + 4 \cdot 4 + 4 \cdot 1) = 142\,020$ parameters (all the weights of all the $45$ BNNs) instead of the $784 \cdot 16 + 16 \cdot 16 + 16 \cdot 10 = 12\,960$ parameters used in \citet{Icarte} for a single large BNN. Note that, in this case, the dimension of the parameter space is $3^{12\,960} (\cong 10^{6\,183}) $, while, in our case, it is $45 \cdot 3^{3\,156} (\cong 10^{1\,507} )$. In the first case, the solver has to find an optimal solution between all the $10^{6\,183}$ different parameter configurations, while with the \emph{BeMi} ensemble, the solver has to find $45$ optimal solutions, each of which lives in a set of cardinality $10^{1\,507}$. This significantly improves the solver performances. We remark that a parameter configuration is given by a weight assignment $\hat{W} = (\hat{w}_{ijl})_{ijl}$ since every other variable is uniquely determined by $\hat{W}$.

	Figure~\ref{fig:comparison} compares the results of our \emph{BeMi} ensemble with the four other methods presented above. Note that 
	the pure MIP model in \citet{Icarte} can handle a single image per class in the given time limit, and so only one point is reported, and note also that
	for the minimum hinge model {\tt MH} presented in \citet{Yorke} only the experiment with 10 digits per class was performed.
	We report the best results reported in the original papers for these four methods.
	
	The \emph{BeMi} ensemble obtains an average accuracy of 68\%, %prima 61%
	outperforms all other approaches when 2, 6 or 10 digits per class are used. 
	Note that our method attains $100\%$ accuracy on the training set, that is, the \texttt{SM} model correctly classifies all the images. In this case, the first model is not needed to ensure feasibility, but it serves mainly as a warm start for the \texttt{MM} model.
	
	\subsection{Experiment 3}
	The goal of the third experiment is to study how our approach scales with the number of data points (i.e., images) per class, and how it is affected by the architecture of the small BNNs within the \emph{BeMi} ensemble.
	For the number of data points per class, we use $10, 20, 30, 40$ training images per digit. 
	We use the layers $\mathcal{N}_a=[784, 4, 4, 1]$ and $\mathcal{N}_b=[784, 10, 3, 1]$ for the two architectures.
	While the first architecture is chosen as to be consistent with the previous experiments, the second one can be described as an integer approximation of $[N, \log_2 N, \log_2 (\log_2 N), \log_2 (\log_2 (\log_2 N))]$.
	Herein, we refer to Experiments $3a$ and $3b$ as the two subsets of experiments related to the architectures $\mathcal{N}_a$ and $\mathcal{N}_b$.
	In both cases, we train each of our $45$ BNNs with a time limit of 290s for model {\tt SM}, 290s for {\tt MM}, and 20s for {\tt MW}, for a total of 600s (i.e., 10 minutes for each BNN).
	
\begin{figure}[t!]
		\centering
		\subfloat[t!][MNIST.\label{fig:mnist}]{%
			%%%%%%%%%%%%%% NEWPICTURE
			
			\begin{tikzpicture}
				\begin{axis}[
					width=6.5cm,
					height=4cm,
					scale only axis,
					xmin=8, xmax=42,
					xtick={10, 20, 30, 40},
					xticklabels={\small 10,\small 20,\small 30, \small 40},
					xmajorgrids,
					ymin=60, ymax=85,
					ytick={60, 65, 70, 75, 80, 85},
					yticklabels={\small 60\%, \small 65\%, \small 70\%, \small 75\%, \small 80\%, \small 85\%},
					ylabel={\small Test Accuracy},
					xlabel={\small Training Images per Digit},
					ymajorgrids,
					axis lines*=left,
					line width=1.0pt,
					mark size=3.0pt,
					legend style ={ at={(1.03,1)}, 
						anchor=north west, %draw=black, 
						fill=white,align=left},
					cycle list name=black white,
					smooth
					]
					
					\fill[mycol1, opacity = 0.1] (10, 67.89) -- (20, 73.35) -- (30, 80.10) -- (40, 79.68) -- (40, 82.54) -- (30, 81.38) -- (20, 77.59) -- (10, 74.54);
					\fill[mycol2, opacity = 0.1] (10, 68.53) -- (20, 73.54) -- (30, 80.30) -- (40, 80.19) -- (40, 82.55) -- (30, 81.41) -- (20, 78.09) -- (10, 74.73);
					
					\addplot[color=mycol1, mark=o, dotted, mark options={solid}] coordinates{
						(10, 70.12)
						(20, 75.37)
						(30, 80.90)
						(40, 81.66)
					};
					\addplot[color=mycol2, mark=square, dashed, mark options={solid}] coordinates{
						(10, 70.36)
						(20, 75.67)
						(30, 80.97)
						(40, 81.75)
					};

				\end{axis}
			\end{tikzpicture}%
		}
		\hfill
		\subfloat[t!][Fashion-MNIST.\label{fig:fashion}]{%
			%%%%%%%%%%%%%%%%%%% NEW PICTURE
			\begin{tikzpicture}
				\begin{axis}[
					width=6.5cm,
					height=4cm,
					scale only axis,
					xmin=8, xmax=42,
					xtick={10, 20, 30, 40},
					xticklabels={\small 10,\small 20,\small 30, \small 40},
					xmajorgrids,
					ymin=60, ymax=85,
					ytick={60, 65, 70, 75, 80, 85},
					yticklabels={\small 60\%, \small 65\%, \small 70\%, \small 75\%, \small 80\%, \small 85\%},
					ylabel={\small Test Accuracy},
					xlabel={\small Training Images per Digit},
					ymajorgrids,
					axis lines*=left,
					line width=1.0pt,
					mark size=3.0pt,
					legend style ={ at={(34/36,1)}, 
						anchor=north east, %draw=black, 
						fill=white,align=left},
					cycle list name=black white,
					smooth
					]
					
					\fill[mycol1, opacity = 0.1] (10, 64.08) -- (20, 69.96) -- (30, 69.13) -- (40, 68.61) -- (40, 72.31) -- (30, 71.29) -- (20, 71.13) -- (10, 71.10);
					\fill[mycol2, opacity = 0.1] (10, 64.24) -- (20, 69.81) -- (30, 68.63) -- (40, 66.73) -- (40, 71.39) -- (30, 71.84) -- (20, 71.46) -- (10, 71.13);
					
					\addplot[color=mycol1, mark=o, dotted, mark options={solid}] coordinates{
						(10, 66.89)
						(20, 70.74)
						(30, 70.35)
						(40, 70.65)
					};
					\addlegendentry{\small $\mathcal{N}_a$};
					\addplot[color=mycol2, mark=square, dashed, mark options={solid}] coordinates{
						(10, 67.07)
						(20, 70.43)
						(30, 70.52)
						(40, 69.06)
					};
					\addlegendentry{\small $\mathcal{N}_b$};

				\end{axis}
			\end{tikzpicture}%
		}
		\caption{Average accuracy for the \emph{BeMi} ensemble tested on two architectures, namely $\mathcal{N}_a=[784,4,4,1]$ and $\mathcal{N}_b=[784,10,3,1]$, using $10, 20, 30, 40$ images per class.}
		\label{fig:fig}
	\end{figure}
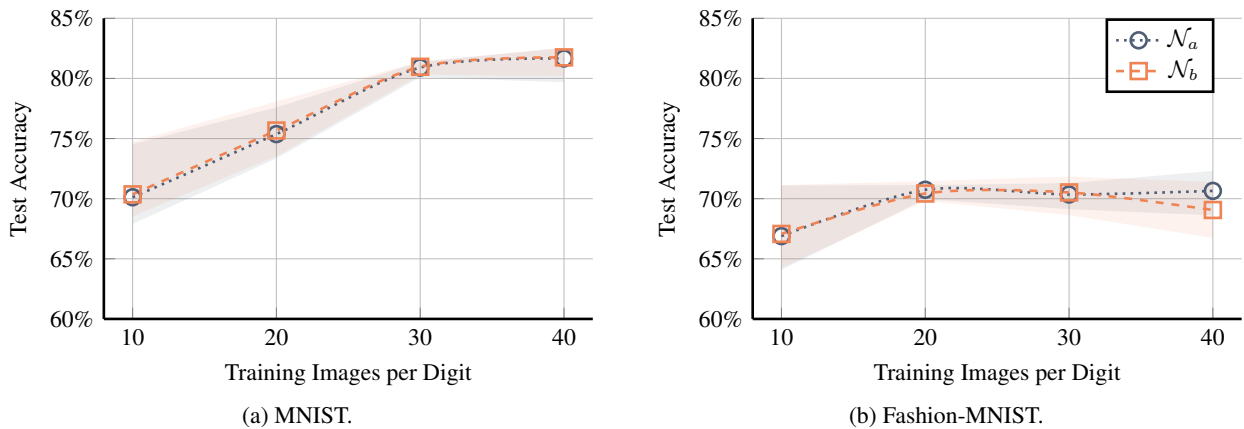

	Figure \ref{fig:mnist} shows the results for Experiments $3a$ and $3b$:
	the dotted and dashed lines refer to the two average accuracies of the two architectures,
	while the coloured areas include all the accuracy values obtained as the training instances vary. 
	The two architectures behave similarly and the best average accuracy exceeds 81\%.
	
	Table~\ref{tab:statustable} reports the results for the \emph{BeMi} ensemble where we distinguish among images classified as correct, wrong, or unclassified. 
	These three conditions refer to different label statuses specified in Definition \ref{def:labels}:
	%Fbemui
	the correct labels are the sum of the statuses $(1C)$ and $(2C)$; the wrong labels of statuses $(2I'), (2I''),$ and $(1I)$; the unclassified labels ({\it n.l.}) of $(oI')$ and $(oI'')$. Notice that the vast majority of the test images have only one dominant label, and so falls into statuses $(1C)$ or $(1I)$. The unclassified images are less than 2.31\%.

	%%%% MNIST AND FASHION-MNIST
	
	\begin{table}[t!]
		\caption{Percentages of MNIST images classified as correct, wrong, or unclassified ({\it n.l.}), and of label statuses, for the architecture $\mathcal{N}_a=[784, 4,4, 1]$. The vast majority of the test images have only one dominant label, and so falls into statuses $(1C)$ or $(1I)$. The unclassified images are less than 2.31\%.}
		\label{tab:statustable}
		\centering
		\setlength{\tabcolsep}{5.4pt}
		\small
		
		\begin{tabular}{c@{\hskip 0.4in}ccc@{\hskip 0.4in}ccccccc}
			\toprule
			Images & \multicolumn{3}{c@{\hskip 0.4in}}{Classification (\%)} & \multicolumn{7}{c}{Label status (\%)} \\
			per class & correct & wrong & {\it n.l.} & $(1C)$ & $(1I)$ & $(2C)$ & $(2I')$ & $(2I'')$ & $(oI')$ & $(oI'')$ \\ \midrule
			10     & 70.12   & 27.65  & 2.23  & 68.43   & 24.53    & 1.69    & 1.21    & 1.91    & 1.88    & 0.35  \\
			20     & 75.37   & 22.32  & 2.31  & 73.79   & 19.33    & 1.58    & 1.39    & 1.60    & 2.02    & 0.29  \\
			30     & 80.90   & 17.46  & 1.64  & 79.64   & 15.01    & 1.26    & 1.25    & 1.20    & 1.46    & 0.18  \\
			40     & 81.66   & 16.68  & 1.66  & 80.34   & 14.36    & 1.32    & 1.09    & 1.23    & 1.45    & 0.21 \\
			\bottomrule
		\end{tabular}
	\end{table}

	%%%%%%%%%%%%%%%%%%%%%%%%%%%%%%%%%%%%%%%%%%%%%%%%
	\subsection{Experiment 4}
	
	\begin{table}[t!]
		\caption{Aggregate results for Experiments 2 and 3: the $4$-th column reports the runtime to solve the first model {\tt SM}; {\it Gap (\%)} refers to the mean and maximum percentage gap at the second MILP model {\tt MM}; {\it Links (\%)} is the percentage of non-zero weights after the solution of models {\tt MM} and {\tt MW}; {\it Active links} is the total number of non-zero weights after the solution of model {\tt MW}, as always averaged over five instances.}
		% \smallskip
		\centering 
		\setlength{\tabcolsep}{5.4pt}
		\small
		\begin{tabular}{c@{\hskip 0.3in}c@{\hskip 0.3in}c@{\hskip 0.3in}c@{\hskip 0.3in}r@{\hskip 0.3in}rr@{\hskip 0.3in}rr@{\hskip 0.2in}r}
			\toprule
			\multirow{2}{*}{Dataset} & \multirow{2}{*}{Layers} & Total & Images & \multicolumn{1}{c@{\hskip 0.3in}}{Model {\tt SM}} & \multicolumn{2}{c@{\hskip 0.3in}}{Gap (\%) } & \multicolumn{2}{c@{\hskip 0.2in}}{Links (\%)} & \multicolumn{1}{c}{Active} \\ 
			& & links & per class & \multicolumn{1}{c@{\hskip 0.3in}}{time (s)} & mean & max & {\tt (MM)} & {\tt (MW)} & \multicolumn{1}{c}{links}\\
			\midrule
			\multirow{8}{*}{MNIST} & \multirow{4}{*}{784,4,4,1} &\multirow{4}{*}{3156} &  10 &  3.00 & 12.06 & 20.70 & 49.13 & 29.21 & 921.80\\
			& & &  20  & 6.47 & 14.81 & 22.18 & 54.70 & 28.41 & 896.60\\
			& & &  30  & 10.60 & 16.04 & 24.08 & 56.44 & 30.46 & 961.40\\
			& & & 40  & 15.04 & 15.98 & 26.22 & 57.92 & 29.27 & 923.80\\[2mm]
			%\cline{2-9}
			& \multirow{4}{*}{784,10,3,1} & \multirow{4}{*}{7873} & 10 & 6.04 & 4.52 & 7.37 & 49.28 & 24.72 & 1946.20\\
			& &  & 20 & 15.01 & 5.46 & 8.40 & 54.97 & 24.40 & 1921.00 \\
			& &  & 30 & 22.68 & 5.92 & 11.00 & 56.73 & 26.96 & 2122.60\\
			& &  & 40 & 33.97 & 6.21 & 20.87 & 58.29 & 24.66 & 1941.40 \\
			\midrule
			\multirow{8}{*}{F-MNIST} & 
			\multirow{4}{*}{784,4,4,1} & \multirow{4}{*}{3156} & 10  & 4.83 & 13.34 & 26.48 & 87.75 & 58.76 & 1854.40 \\
			& &  & 20  & 9.77 & 14.05 & 28.91 & 90.73 & 59.97 & 1892.60 \\
			& & &  30  & 36.10 & 19.95 & 136.33 & 92.41 & 58.12 & 1834.20 \\
			& & & 40  & 72.15 & 30.36 & 333.70 & 93.57 & 59.46 & 1876.60 \\[2mm]
			%\cline{2-9}  
			& \multirow{4}{*}{784,10,3,1} & \multirow{4}{*}{7873}&  10 & 11.42 & 4.87 & 9.14 &  87.90 & 51.57 & 4060.20 \\
			& & &  20 &  21.46 & 5.11 & 9.86 & 91.05 & 52.29 & 4116.80\\
			& & &  30 &  35.12 & 6.30 & 40.07 & 92.78 & 52.62 & 4142.80\\
			& & & 40 &  52.37 & 8.99 & 56.14 & 94.01 & 53.38 & 4202.60\\
			\bottomrule  
		\end{tabular}
		\label{tab:all_stages}
		\vspace{-10pt}
	\end{table}
	
	The goal of the fourth experiment is to revisit the questions of Experiments~$3a$ and $3b$ with the two architectures $\mathcal{N}_a$ and $\mathcal{N}_b$, using the more challenging Fashion-MNIST dataset.

	Figure~\ref{fig:fashion} shows the results of Experiments $3a$ and $3b$. 
	As in Figure \ref{fig:comparison}, the dotted and dashed lines represent the average percentages of correctly classified images, while the coloured areas include all accuracy values obtained as the instances vary.
	For the Fashion-MNIST, the best average accuracy exceeds 70\%.
	
	Table \ref{tab:all_stages} reports detailed results for all Experiments 2 and 3. The first two columns give the dataset and the architecture, the third column reports the total number of links, i.e., the total number $w$ variables, for each architecture, and the fourth column specifies the number of images per digit used during training.
	The $5$-th column reports the runtime for solving model {\tt SM}.
	Note that the time limit is 290 seconds; hence, we solve exactly the first model, consistently  achieving a training accuracy of 100\%.
	The remaining five columns give:
	{\it Gap (\%)} refers to the mean and maximum percentage gap at the second MILP model {\tt (MM)}
	of our lexicographic multi-objective model, as reported by the Gurobi {\tt MIPgap} attribute;
	{\it Links (\%)} is the percentage of non-zero weights after the solution of the second model {\tt MM}, and after the solution of the last model {\tt MW}.
	{\it Active links} is the total number of non-zero weights after the solution of the last model {\tt MW}.
	The results show that the runtime and the gap increase with the size of the input set. 
	However, for the percentage of removed links, there is a significant difference between the two datasets: for MNIST, our third model {\tt MW} removes around 70\% of the links, while for the Fashion-MNIST, it removes around 50\% of the links.
	Note that in both cases, these significant reductions show how our model is also optimizing the BNN architecture.
	Furthermore, note that even if the accuracy of the two architectures is comparable, the total number of non-zero weights of $\mathcal{N}_a$ is about half the number of non-zero weights of $\mathcal{N}_b$.
	
	\subsection{Experiment 5}
	
	\begin{figure}[t!]
		\centering
		\begin{tikzpicture}
			\begin{axis}[
				width=14cm,
				height=6cm,
				scale only axis,
				xmin=0, xmax=42,
				xtick={2,6,10, 20, 30, 40},
				xticklabels={\small 2,\small 6,\small 10, \small 20, \small 30, \small 40},
				xmajorgrids,
				ymin=45, ymax=95,
				ytick={45, 55, 65, 75, 85, 95},
				yticklabels={\small 45\%, \small 55\%, \small 65\%, \small 75\%, \small 85\%, \small 95\%,},
				ylabel={\small Test Accuracy},
				xlabel={\small Training Images per Digit},
				legend cell align={right},
				ymajorgrids,
				axis lines*=left,
				line width=1.0pt,
				mark size=3.0pt,
				legend style ={ at={(40/42,0.6)}, 
					anchor=north east, %draw=black, 
					fill=white,align=left},
				cycle list name=black white,
				smooth
				]

				\fill[mycol1, opacity = 0.1] (2, 48.56) -- (6, 66.13) -- (10, 62.19) -- (20, 79.81) -- (30, 83.63) -- (40, 85.00) -- (40, 90.75) -- (30, 89.50) -- (20, 87.19) -- (10, 84.94) -- (6, 76.44) -- (2, 68.69);
				
				\fill[mycol2, opacity = 0.1] (2, 48.25) -- (6, 66.69) -- (10, 61.44) -- (20, 79.88) -- (30, 83.75) -- (40, 85.44) -- (40, 90.81) -- (30, 89.44) -- (20, 87.75) -- (10, 84.38) -- (6, 75.88) -- (2, 68.31);
				
				\fill[mycol3, opacity = 0.1] (2, 48.94) -- (6, 66.25) -- (10, 62.19) -- (20, 80.56) -- (30, 84.38) -- (40, 85.69) -- (40, 90.63) -- (30, 90.13) -- (20, 87.19) -- (10, 84.63) -- (6, 76.56) -- (2, 69.09);
				
				\fill[mycol4, opacity = 0.1] (2, 48.56) -- (6, 66.75) -- (10, 61.69) -- (20, 80.06) -- (30, 84.19) -- (40, 85.44) -- (40, 90.69) -- (30, 89.81) -- (20, 87.19) -- (10, 84.31) -- (6, 76.81) -- (2, 68.94);

				\addplot[color=mycol1, mark=o, dotted, mark options={solid}] coordinates{
					(2, 60.43)
					(6, 72.48)
					(10, 71.45)
					(20, 82.79) 
					(30, 87.11) 
					(40, 87.90) 
				};
				\addlegendentry{\small $w_{ilj} \in \{-1, 0, 1\} $};
				\addplot[color=mycol2, mark=square, dashed, mark options={solid}] coordinates{
					(2, 60.25)
					(6, 72.14)
					(10, 71.28)
					(20, 82.89) 
					(30, 86.89) 
					(40, 87.86) 
				};
				\addlegendentry{\small $w_{ilj} \in \{-3, \dots, 3\}$};
				\addplot[color=mycol3, mark=triangle, dash dot, mark options={solid}] coordinates{
					(2, 60.53)
					(6, 72.21)
					(10, 71.45)
					(20, 83.05) 
					(30, 87.16) 
					(40, 87.75) 
				};
				\addlegendentry{\small $w_{ilj} \in \{-7, \dots, 7\}$};
				\addplot[color=mycol4, mark=star, densely dotted, mark options={solid}] coordinates{
					(2, 60.58)
					(6, 72.23)
					(10, 71.35)
					(20, 83.15) 
					(30, 87.14) 
					(40, 87.75) 
				};
				\addlegendentry{\small $w_{ilj} \in \{-15, \dots, 15\}$};
				
			\end{axis}
		\end{tikzpicture}%

		\caption{Comparison of accuracy for different values of $P$. Note how using an exponentially larger research space, namely, using values of $P$ greater than $1$, does not improve the accuracy.}    
		\label{fig:bnn_vs_inn}
	\end{figure}
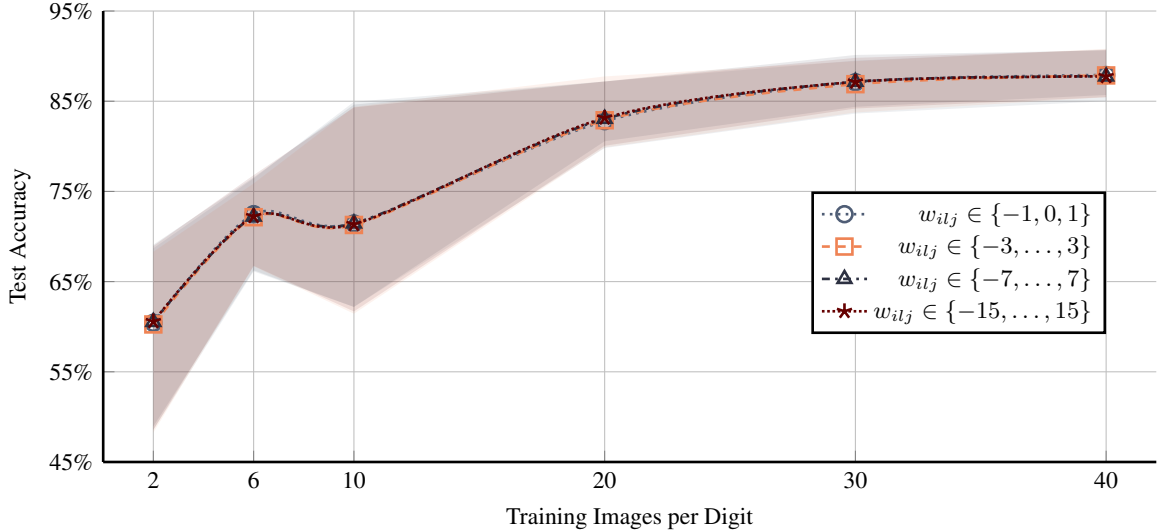
	\begin{table}[t!]
		\caption{Weights distributions of the INNs whose accuracy is depicted in Figure \ref{fig:bnn_vs_inn}. The column $w=-P$ indicates the percentage of weights that are equal to $-P$, and so on. The networks have different values but similar extremal weights distributions, where less than $2\%$ of the weights attain a value that is not $P$, $-P$, or zero, indicated as \emph{others}.}
		% \smallskip
		\centering 
		\setlength{\tabcolsep}{5.4pt}
		\small
		\begin{tabular}{c@{\hskip 0.4in}r@{\hskip 0.4in}r@{\hskip 0.4in}r@{\hskip 0.4in}r@{\hskip 0.4in}r}
			\toprule
			value of $P$ & Images per class & $w = -P$ & $w = 0$ & $w = P$ & \emph{others}\\
			\midrule
			\multirow{6}{*}{1} & 2 & 4.06 & 73.66 & 22.28 & - \\
			& 6 & 4.75 & 73.70 & 21.55 & - \\
			& 10 & 5.15 & 74.02 & 20.83 & - \\
			& 20 & 7.12 & 61.31 & 31.57 & - \\
			& 30 & 6.21 & 74.89 & 18.90 & - \\
			& 40 & 13.16 & 64.98 & 21.86 & - \\
			\midrule%\cline{2-6}
			\multirow{6}{*}{3} & 2 & 3.86 & 75.91 & 20.15 & 0.08 \\
			& 6 & 4.60 & 73.46 & 21.63 & 0.31 \\
			& 10 & 6.25 & 73.50 & 19.72 & 0.53 \\
			& 20 & 12.48 & 69.58 & 17.19 & 0.75 \\
			& 30 & 13.06 & 74.00 & 11.89 & 1.05 \\
			& 40 & 13.99 & 65.51 & 19.31 & 1.19 \\
			\midrule%\cline{2-6}
			\multirow{6}{*}{7} & 2 & 3.90 & 75.89 & 20.10 & 0.11 \\
			& 6 & 4.46 & 73.40 & 21.76 & 0.38 \\
			& 10 & 8.82 & 68.00 & 22.60 & 0.58 \\
			& 20 & 12.46 & 65.56 & 21.06 & 0.92 \\
			& 30 & 14.99 & 73.82 & 9.95 & 1.24 \\
			& 40 & 18.05 & 70.90 & 9.40 & 1.65 \\
			\midrule%\cline{2-6}
			\multirow{6}{*}{15} & 2 & 3.69 & 75.89 & 20.30 & 0.12 \\
			& 6 & 4.68 & 73.31 & 21.64 & 0.37 \\
			& 10 & 5.53 & 73.28 & 20.61 & 0.58 \\
			& 20 & 7.43 & 69.34 & 22.19 & 1.04 \\
			& 30 & 13.76 & 67.64 & 17.23 & 1.37 \\
			& 40 & 17.34 & 70.86 & 10.13 & 1.67 \\
			\bottomrule
		\end{tabular}
		\label{tab:weights_dist}
	\end{table}
	
	The goal of the fifth experiment is to compare the performances of BNNs and non-trivial INNs. The results refer to five different runs of an INN of architecture $[784, 4, 4, 1]$ specialised in distinguishing between digits 4 and 9 of the MNIST dataset. The number of training images varies between $2$, $6$, $10$, $20$, $30$, and $40$, while the test images are $800$ per digit, so $1600$ in total. The imposed time limit is $290$s + $290$s + $20$s.
	
	As Figure~\ref{fig:bnn_vs_inn} shows, different INNs are comparable not only in the average accuracy, represented by the dotted lines, but also in the maximum and minimum accuracy, reported by the shaded areas.
	
	In order to study why different values of $P$ lead to comparable accuracy, we report the weights distributions of the INNs whose accuracy is depicted in Figure \ref{fig:bnn_vs_inn}. Table \ref{tab:weights_dist} highlights not only the role of the \texttt{MW} model but also the extreme-valued nature of the distributions. In fact, it can be seen that apart from the percentage of weights set to zero, the vast majority of the remaining weights have a value of either $P$ or $-P$, with less than $1.67\%$ of the weights attaining one of the other intermediate values. We remark that this type of phenomenon is recurrent in the literature under the name of magnitude pruning, see \citet{han2015deep, morcos2019one,blalock2020state}. In our setting, such a phenomenon occurs spontaneously.
	
	\subsection{Experiment 6}
	The goal of the sixth experiment is to study the impact of the multi-objective model and the weight distributions over a different dataset, namely the Heart Disease Dataset by \citet{heart}.
	
	Table \ref{tab:heart} reports the average accuracy and weights distributions of this experiment. The average was performed over 5 instances, and for each instance $200$ data samples were used, where 80\% was used for the training and 20\% was used for the test. All the networks have the same architecture, namely $[13, 5, 1]$, and each network's imposed time limit is $60$ minutes, equally distributed in the steps of each model: $60$ minutes for the \texttt{SM}, $30+30$ for the \texttt{SM}+\texttt{MM}, and $20+20+20$ for \texttt{SM}+\texttt{MM}+\texttt{MW}.
	\begin{table}[t!]
		\caption{Average accuracy and weight distribution for the Heart Disease dataset. The column Average accuracy depicts the average percentage of correctly classified test data. The column $w=-P$ indicates the percentage of weights that are equal to $-P$, and so on. For each value of $P$, the best result in terms of accuracy with respect to the model is written in bold. Notice that in the majority of the cases, the multi-objective function performs better than the single-objective one. Notice also that even if the percentage of the non-zero and non-extremal weights, indicated as \emph{others}, is higher than the one obtained with the MNIST dataset, the distribution is still not uniform.}
		\label{tab:heart}
		% \smallskip
		\centering
		\setlength{\tabcolsep}{5.4pt}
		\small
		\begin{tabular}{cc@{\hskip 0.4in}r@{\hskip 0.4in}rrrr}
			\toprule
			Models & value of $P$ & Average accuracy & $w= -P$ & $w = 0$ & $w=P$ & \emph{others}\\
			\midrule
			\multirow{4}{*}{\texttt{SM}} & 1 & 74.00  & 35.43& 45.14& 19.43&-\\
			& 3 & 75.00   &21.43 & 27.43& 16.57 &34.57\\
			& 7 & \textbf{77.00}  &27.43 &14.86 & 16.86&40.85\\
			& 15 & 77.00   & 18.86& 11.71& 12.86&56.57\\
			\midrule
			\multirow{4}{*}{\texttt{SM$\,+\,$MM}} & 1 & 75.00  & 14.57& 18.29& 67.14&-\\
			& 3 & 75.50  & 10.57& 13.43&63.14 &12.86\\
			& 7 & 73.50   &12.29 & 9.43& 52.57&25.71\\
			& 15 & \textbf{78.50}   & 11.71 & 5.71&46.86&35.72\\
			\midrule
			\multirow{4}{*}{\texttt{SM$\,+\,$MM$\,+\,$MW}} & 1 & \textbf{75.50}  & 9.43& 27.43&63.14&-\\
			& 3 & \textbf{76.00} & 5.43 & 25.43& 57.43&11.71\\
			& 7 & 73.50   &7.43 & 18.57& 49.71&24.29\\
			& 15 & \textbf{78.50}   &8.57 & 13.43& 45.71&32.29\\
			\bottomrule
		\end{tabular}
	\end{table}

	%%%%%%%%%%%%%%%%%%%%%%%%%%%%%%%%%%%%%%%%%%%%%%%
	
	\section{Conclusion and Future Work}
	\label{sec:conc}
	
	This paper introduced the \emph{BeMi} ensemble, a structured architecture of INNs for classification tasks.
	Each network specializes in distinguishing between pairs of classes and combines different approaches already existing in the literature to preserve feasibility while being robust and simple.
	These features and the nature of the parameters are critical to enabling neural networks to run on low-power devices. In particular, binarized NNs can be implemented using Boolean operations and do not need GPUs to run.
	
	The output of the \emph{BeMi} ensemble is chosen by a majority voting system that generalizes the one-versus-one scheme.
	Notice that the \emph{BeMi} ensemble is a general architecture that could be employed using other types of neural networks.
	
	An interesting conclusion from our computational experiments is a counter-intuitive result: that the greater flexibility in the search space
	of INNs does not necessarily result in better classification accuracy compared to BNNs.  We find
	it noteworthy that our computational evidence supports the idea that simpler
	BNNs are either superior or equal to INNs in terms of accuracy.
	
	A current limitation of our approach is the strong dependence on the randomly sampled data used for training.
	In future work, we plan to improve the training data selection by using a $k$-medoids approach, dividing all images of the same class into disjoint non-empty subsets and consider their centroids as training data.
	This approach should mitigate the dependency on the sampled training data points.
	
	Second, we also plan to better investigate the scalability of our method with respect to the number of classes of the classification problem training fewer BNNs, namely, one for every $\mathcal{J} \in \mathcal{Q} \subset \mathcal{P}$, with $|\mathcal{Q}| \ll |\mathcal{P}|$.
	Besides the datasets we exploited, in future, we intend to investigate datasets more appropriate for the task of few-shot learning \citep{brigato}.
	
	Third, another possible future research direction is to exploit the generalisation of our ensemble, allowing to have networks which distinguishes between $m$ classes instead of $2$, where the total number of classes of the problem is $n \gg m$. Note that the definitions of this generalization are presented in Appendix \ref{app:ensemble}.
	
	Fourth, an interesting future research direction regards stochastic/robust optimization and scenario generation, in the following sense. In the case of the MNIST/FashionMNIST for example, the images can be seen as samples of many unknown probability distributions, one for each class: if one would like to train a neural network with a MIP, using few images, the selection of these samples with which the training is performed is crucial, and so the study of this problem from a stochastic point of view could lead to interesting results.
	
	% Acknowledgments here
	\section*{Acknowledgment}
	The authors thank the anonymous reviewers for their comments.  We thank the participants of the Dagstuhl  
	Seminar 22431 `Data-Driven Combinatorial Optimisation', and also thank T\'omas \TH{}orbjarnarson.  This research was partially supported by TAILOR, a project funded by EU Horizon 2020 research and innovation programme under grant number~952215. The work of the author A.M.\ Bernardelli is supported by a Ph.D. scholarship funded under the `Programma Operativo Nazionale Ricerca e Innovazione' 2014--2020.
	
	% Leave this (end of acknowledgment)

	\clearpage
	
	\appendix
	\section{Generalization of the ensemble}\label{app:ensemble}
	The definition of the ensemble introduced in Section \ref{sub32} can be generalize as follows. Define $\mathcal{P}(\mathcal{I})_m$ as the set of all the subsets of the set $\mathcal{I}$ that have cardinality $m$, where $\mathcal{I}$ is the set of the classes of the classification problem.  Then our structured ensemble is constructed in the following way:
	\begin{enumerate}
		\item We set a parameter $1 < m \leq n = |\mathcal{I}|$.
		
		\item We train a INN denoted by $\mathcal{N}_{\mathcal{J}}$ for every $\mathcal{J} \in \mathcal{P}(\mathcal{I})_{m}$.
		
		\item When testing a data point, we feed it to our list of trained INNs, namely $(\mathcal{N}_{\mathcal{J}})_{\mathcal{J} \in \mathcal{P}(\mathcal{I})_{m}}$, obtaining a list of predicted labels $(\mathfrak{e}_{\mathcal{J}})_{\mathcal{J} \in \mathcal{P}(\mathcal{I})_{m}}$.
		
		\item We then apply a majority voting system.
	\end{enumerate} 
	
	Note that we set $m >1$, otherwise our structured ensemble would have been meaningless. 
	Whenever $m = n$, our ensemble is made of one single INN.
	When $m = 2$, we are using the one-versus-one scheme presented in the paper.
	
	The idea behind this structured ensemble is that, given an input $\bm{x}^k$ labelled $l$ $(=y^k)$, the input is fed into $\binom{n}{m}$ networks where $\binom{n-1}{m-1}$ of them are trained to recognize an input with label $l$. 
	If all of the networks correctly classify the input $\bm{x}^k$, then at most $\binom{n-1}{m-1} - \binom{n-2}{m-2}$ other networks can classify the input with a different label. 
	With this approach, if we plan to use $r \in \mathbb{N}$ inputs for each label, we are feeding our INNs a total of $m \cdot r$ inputs instead of feeding $n \cdot r$ inputs to a single large INN. 
	When $m  \ll n$, it is much easier to train our structured ensemble of INNs rather than training one large INN. 

	\paragraph{Majority voting system.}
	After the training, we feed one input $\bm{x}^k$ to our list of INNs, and we need to elaborate on the set of outputs.
	
	\begin{definition}[Dominant label]
		For every $b \in \mathcal{I}$, we define
		\begin{equation*}
			C_b = \{\mathcal{J} \in \mathcal{P}(\mathcal{I})_m \; | \;  \mathfrak{e}_{\mathcal{J}} = b\},
		\end{equation*}
		and we say that a label $b$ is a \emph{dominant label} if  $|C_b| \geq |C_l|$ for every $l \in \mathcal{I}$. 
		We then define the set of dominant labels
		\begin{equation*}
			\mathcal{D} \coloneqq \{b \in \mathcal{I} \; | \; b \text{ is a dominant label}\}.
		\end{equation*}
	\end{definition}
	
	Using this definition, we can have three possible outcomes:
	\begin{itemize}
		\item[(a)] There exists a label $b \in \mathcal{I}$ such that $\mathcal{D}= \{b\}$ $\implies$ our input is labelled as $b$.
		\item[(b)] There exist $b_1, \dots, b_m \in \mathcal{I}$, $b_i \neq b_j$ for all $i \neq j$ such that $\mathcal{D} = \{b_1, \dots, b_m\}$, so $\mathcal{D} \in \mathcal{P}(\mathcal{I})_m$  $\implies$ our input is labelled as $\mathfrak{e}_{\{b_1,  \dots, b_m\}} = \mathfrak{e}_{\mathcal{D}}$.  
		\item[(c)] $|\mathcal{D}| \neq 1 \land |\mathcal{D}| \neq  m$ $\implies$ our input is labelled as $z \notin \mathcal{I}$.
	\end{itemize}
	While case (a) is straightforward, we can label our input even when we do not have a clear winner, that is, when we have trained a INN on the set of labels that are the most frequent (i.e., case (b)).
	Note that the proposed structured ensemble alongside its voting scheme can also be exploited for regular NNs.
	\begin{definition}[Label statuses]
		In our labelling system, when testing an input seven different cases, herein called \emph{label statuses}, can arise. The statuses names are of the form \char96\char96 number of the dominant labels + fairness of the prediction". The first parameter can be $1, m$, or $o$, where $o$ means "other cases". The fairness of the prediction is $C$ when it is correct, or $I$ when it is incorrect. The subscripts related to $I'$ and $I''$ only distinguish between different cases. These cases are described through the following tree diagram. \vspace*{15pt}
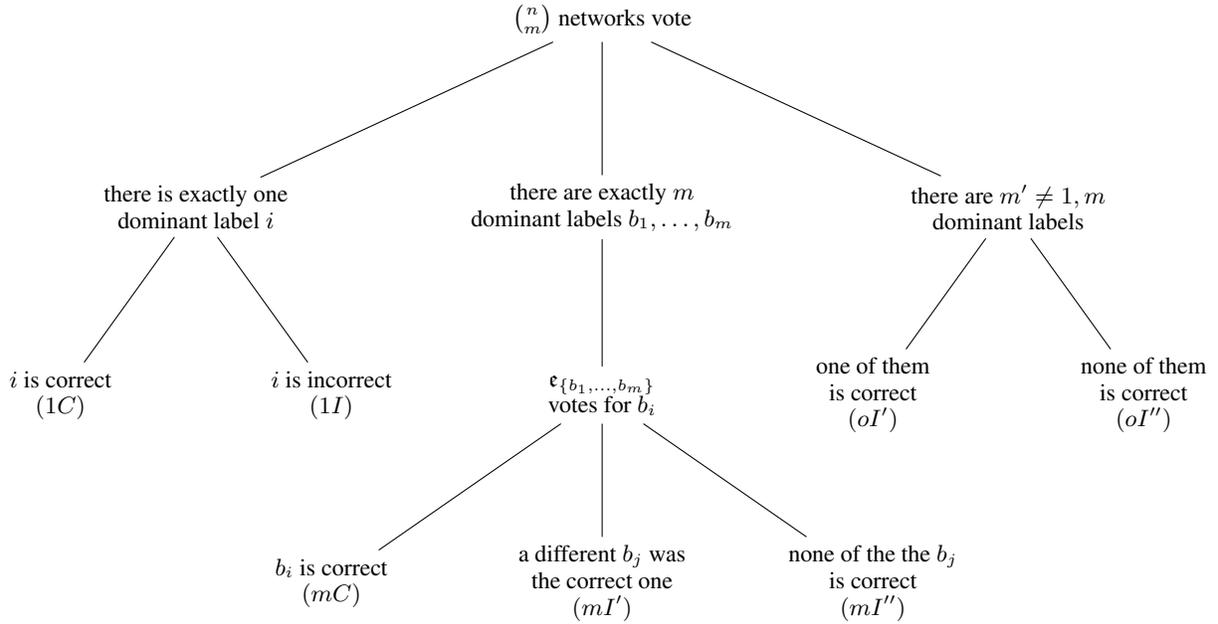
		\begin{center}
			\small
			\begin{tikzpicture}[level distance=2.5cm,
				level 1/.style={sibling distance=5.4cm},
				level 2/.style={sibling distance=3.6cm}]
				\node { $\binom{n}{m}$ networks vote}
				child {node[align=center] {there is exactly one\\ dominant label $i$}
					child {node[align=center] {$i$ is correct \\ $(1C)$ }}
					child {node[align=center] {$i$ is incorrect \\ $(1I)$}}}
				child {node[align=center] {there are exactly $m$ \\ dominant labels $b_1, \dots, b_m$}
					child {node[align=center] {$ \mathfrak{e}_{\{b_1, \dots, b_m\}}$ \\ votes for $b_i$}
						child {node[align=center] {$b_i$ is correct \\ $(mC)$}}
						child {node[align=center] {a different $b_j$ was\\  the correct one \\ $(mI')$}}
						child {node[align=center] {none of the the $b_j$\\ is correct\\ $(mI'')$}}}}
				child {node[align=center] {there are $m' \neq 1, m$ \\ dominant labels}
					child {node[align=center] {one of them \\is correct \\ \textbf{$(oI')$}}}
					child {node[align=center] {none of them \\is correct \\ \textbf{$(oI'')$}}}};
			\end{tikzpicture}
		\end{center}
		
	\end{definition}
	The reason behind this generalisation is the following. Suppose to have trained a NN to distinguish between $3$ different classes, $c_1, c_2, c_3$. If one class is added to the problem, namely $c_4$, instead of discarding the trained NN, one could train $3$ other NNs, namely the one that distinguishes between $c_1, c_2, c_4$, another one for $c_1, c_3, c_4$, and the last one for $c_2, c_3, c_4$, and use the majority voting scheme. Note that the training of the smaller networks is linked to smaller MILPs: in fact, if we plan to use $10$ input data for each class, we need $40$ input data in total, but every network is fed with only $30$ of them.  Moreover, the training of the three additional NNs can be done in parallel, saving computational time.

	\section{Complement to Experiment 1}
	\label{app:exp}
	
	\begin{figure}[t!]
		\centering
		\begin{tikzpicture}[scale=0.9]
			\begin{axis}[height=6.7cm, width=0.9\textwidth,
				colormap={orangewhite}{color=(white);color=(mycol2!80)},
				xlabel=Predicted,
				xlabel style={yshift=-3pt},
				ylabel=Actual,
				ylabel style={yshift=10pt},
				xticklabels={$0$, $1$, $2$, $3$, $4$, $5$, $6$, $7$, $8$, $9$, $n.l.$}, % changed
				xtick={0,1,2,3,4,5,6,7,8,9,10}, % changed
				xtick style={draw=none},
				yticklabels={$0$, $1$, $2$, $3$, $4$, $5$, $6$, $7$, $8$, $9$}, % changed
				ytick={0,1,2,3,4,5,6,7,8,9}, % changed
				ytick style={draw=none},
				enlargelimits=false,
				colorbar,
				xticklabel style={
					rotate=0
				},
				nodes near coords={\pgfmathprintnumber\pgfplotspointmeta},
				nodes near coords style={
					yshift=-7pt
				},
				]
				\addplot[
				matrix plot,
				mesh/cols=11, % changed
				point meta=explicit,draw=gray,
				table/row sep=\\]
				table[meta=C]{
					x y C\\
					0 0 89.6\\
					1 0 0\\
					2 0 0\\
					3 0 0.4\\
					4 0 0.2\\
					5 0 3.4\\
					6 0 1.2\\
					7 0 0\\
					8 0 1\\
					9 0 0.2\\
					10 0 4\\
					0 1 0\\
					1 1 97.4\\
					2 1 0.4\\
					3 1 0.8\\
					4 1 0\\
					5 1 0.2\\
					6 1 0.4\\
					7 1 0\\
					8 1 0.4\\
					9 1 0.2\\
					10 1 0.2\\ 
					0 2 2.4\\
					1 2 11.6\\
					2 2 67.8\\
					3 2 3.8\\
					4 2 0.4\\
					5 2 1.6\\
					6 2 3.4\\
					7 2 0.4\\
					8 2 5\\
					9 2 0.4\\
					10 2 3.2\\
					0 3 0\\
					1 3 1.8\\
					2 3 2.2\\
					3 3 81.2\\
					4 3 0\\
					5 3 3.8\\
					6 3 0.4\\
					7 3 2.4\\
					8 3 4.8\\
					9 3 1\\
					10 3 2.4\\
					0 4 0\\
					1 4 0.2\\
					2 4 2.8\\
					3 4 0\\
					4 4 76.4\\
					5 4 0.2\\
					6 4 3.8\\
					7 4 1.2\\
					8 4 0.8\\
					9 4 13.2\\
					10 4 1.4\\
					0 5 0.6\\
					1 5 1\\
					2 5 0.2\\
					3 5 7\\
					4 5 0.4\\
					5 5 80\\
					6 5 1.4\\
					7 5 0.6\\
					8 5 5\\
					9 5 1.6\\
					10 5 1.2\\
					0 6 2.4\\
					1 6 0.6\\
					2 6 3\\
					3 6 0.2\\
					4 6 1.8\\
					5 6 4.6\\
					6 6 83.6\\
					7 6 0.4\\
					8 6 0.6\\
					9 6 0\\
					10 6 1.8\\
					0 7 0.2\\
					1 7 4\\
					2 7 5.2\\
					3 7 2.8\\
					4 7 1.4\\
					5 7 0.2\\
					6 7 0.2\\
					7 7 76.6\\
					8 7 1.2\\
					9 7 7\\
					10 7 1.2\\
					0 8 1.6\\
					1 8 2.6\\
					2 8 1.8\\
					3 8 6.2\\
					4 8 1.4\\
					5 8 6.4\\
					6 8 0.8\\
					7 8 1\\
					8 8 73.4\\
					9 8 2.4\\
					10 8 2.4\\
					0 9 1.4\\
					1 9 1\\
					2 9 1\\
					3 9 1.2\\
					4 9 8.6\\
					5 9 0\\
					6 9 0\\
					7 9 3\\
					8 9 1.6\\
					9 9 80.4\\
					10 9 1.8\\  
				}; % added every entry where x=4 or y=4
			\end{axis}
			
		\end{tikzpicture}
		\caption{Confusion matrix of the BeMi ensemble trained with $40$ images per digits with the architecture $[784, 4, 4, 1]$. We reported the percentages for each couple. Note that the $n.l.$ indicate the unclassified images.}    
		\label{fig:confusion_matrix}
	\end{figure}
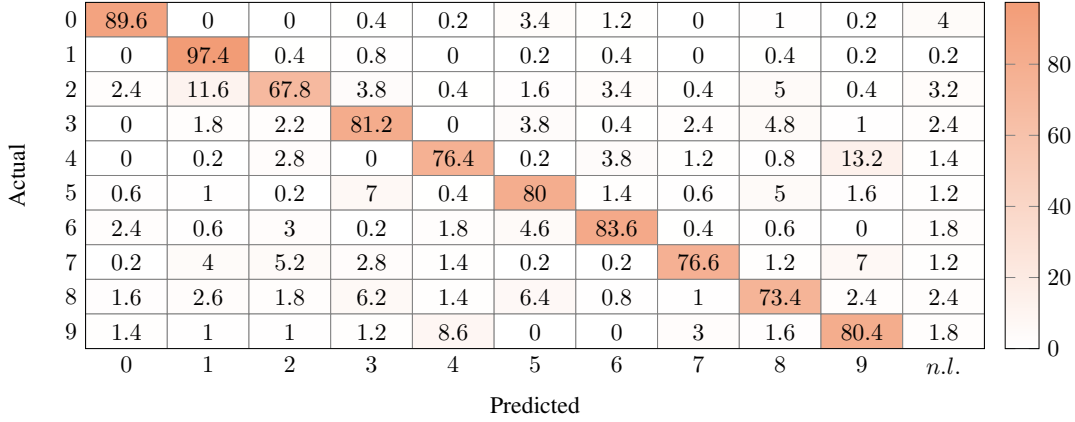
	
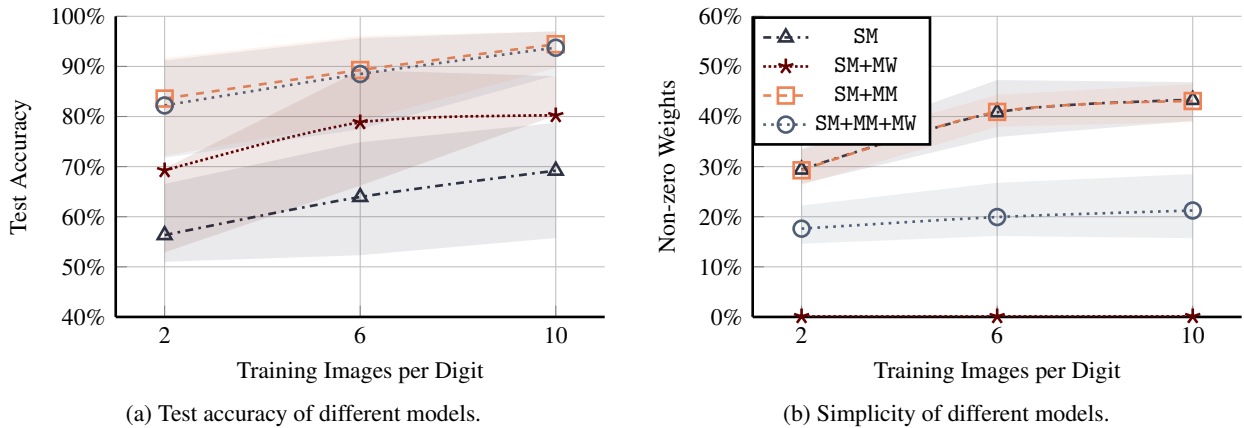
\begin{figure}[t!]
		\centering
		\subfloat[t!][Test accuracy of different models.\label{fig:acc_abl18}]{%
			\begin{tikzpicture}
				
				\begin{axis}[
					width=6.5cm,
					height=4cm,
					scale only axis,
					xmin=1, xmax=11,
					xtick={2,6,10},
					xticklabels={\small 2,\small 6,\small 10},
					xmajorgrids,
					ymin=40, ymax=100,
					ytick={40, 50, 60, 70, 80, 90, 100},
					yticklabels={\small 40\%, \small 50\%, \small 60\%, \small 70\%, \small 80\%, \small 90\%, \small 100\%},
					ylabel={\small Test Accuracy},
					xlabel={\small Training Images per Digit},
					legend cell align={left},
					ymajorgrids,
					axis lines*=left,
					line width=1.0pt,
					mark size=3.0pt,
					legend style ={ at={(1,1)},
						anchor=north west, 
						fill=white,align=left},
					cycle list name=black white,
					smooth
					]
					
					\addplot[color=mycol3, mark=triangle, dash dot, mark options={solid}] coordinates{
						(2, 56.33)
						(6, 63.94)
						(10, 69.19) 
					};

					\addplot[color=mycol4, mark=star, densely dotted, mark options={solid}] coordinates{
						(2, 69.23)
						(6, 78.81)
						(10, 80.21)
					};

					\addplot[color=mycol2, mark=square, dashed, mark options={solid}] coordinates{
						(2, 83.58)
						(6, 89.28)
						(10, 94.43) 
					};

					\addplot[color=mycol1, mark=o,dotted, mark options={solid}] coordinates{
						(2, 82.21)
						(6, 88.48)
						(10, 93.73) 
					};
					
					\fill[mycol3, opacity = 0.1] (2, 51) -- (6, 52.31) -- (10, 55.75) -- (10, 78.75) -- (6, 74.88) -- (2, 66.56);
					\fill[mycol4, opacity = 0.1] (2, 52.88) -- (6, 66.19) -- (10, 80.21) -- (10, 88.06) -- (6, 89.25) -- (2, 69.65);
					\fill[mycol2, opacity = 0.1] (2, 72.13) -- (6, 79.06) -- (10, 89.69) -- (10, 97.00) -- (6, 96.06) -- (2, 91.63);
					\fill[mycol1, opacity = 0.1] (2, 71.75) -- (6, 77.43) -- (10, 88.25) -- (10, 97.00) -- (6, 95.63) -- (2, 91.19);
					
				\end{axis}
			\end{tikzpicture}%
		}\hfill
		%%%%%%%%%%%%%% TDOD da fare distriuzione pesi
		\subfloat[t!][Simplicity of different models.\label{fig:weights_abl18}]{%
			\begin{tikzpicture}
				\begin{axis}[
					width=6.5cm,
					height=4cm,
					scale only axis,
					xmin=1, xmax=11,
					xtick={2,6,10},
					xticklabels={\small 2,\small 6,\small 10},
					xmajorgrids,
					ymin=0, ymax=60,
					ytick={0, 10, 20, 30, 40, 50, 60},
					yticklabels={\small 0\%, \small 10\%, \small 20\%, \small 30\%, \small 40\%, \small 50\%, \small 60\%,},
					ylabel={\small Non-zero Weights},
					xlabel={\small Training Images per Digit},
					ymajorgrids,
					axis lines*=left,
					line width=1.0pt,
					mark size=3.0pt,
					legend style ={ at={(0,1)}, 
						anchor=north west, draw=black, 
						fill=white,align=left},
					cycle list name=black white,
					smooth
					]
					\addplot[color=mycol3, mark=triangle, dash dot, mark options={solid}] coordinates{
						(2, 29.41)
						(6, 40.84)
						(10, 43.31) 
					};
					\addlegendentry{\small \texttt{SM}};
					\addplot[color=mycol4, mark=star, densely dotted, mark options={solid}] coordinates{
						(2, 0.1)
						(6, 0.1)
						(10, 0.11)
					};
					\addlegendentry{\small \texttt{SM}+\texttt{MW}};
					\addplot[color=mycol2, mark=square, dashed, mark options={solid}] coordinates{
						(2, 29.27)
						(6, 40.92)
						(10, 43.09) 
					};
					\addlegendentry{\small \texttt{SM}+\texttt{MM}};
					\addplot[color=mycol1, mark=o,dotted, mark options={solid}] coordinates{
						(2, 17.62)
						(6, 19.94)
						(10, 21.24) 
					};
					\addlegendentry{\small \texttt{SM}+\texttt{MM}+\texttt{MW}};
					\fill[mycol3, opacity = 0.1] (2, 26.58) -- (6, 35.90) -- (10, 39.13) -- (10, 46.86) -- (6, 47.28) -- (2, 33.52);
					\fill[mycol4, opacity = 0.1] (2, 0.1) -- (6, 0.1) -- (10, 0.1) -- (10, 0.13) -- (6, 0.13) -- (2, 0.1);
					\fill[mycol2, opacity = 0.1] (2, 26.33) -- (6, 37.96) -- (10, 38.89) -- (10, 46.51) -- (6, 44.33) -- (2, 33.30);
					\fill[mycol1, opacity = 0.1] (2, 14.61) -- (6, 16.16) -- (10, 15.72) -- (10, 28.52) -- (6, 26.77) -- (2, 22.24);
					
				\end{axis}
			\end{tikzpicture}%
		}
		\caption{Ablation study with a different couple of digits, namely $1$ and $8$, which are easier to distinguish.}
		\label{fig:ablation18}
	\end{figure}

	We first performed an experiment with the \emph{BeMi} ensemble trained with $40$ images per digit and we reported the results in a confusion matrix, depicted in Figure \ref{fig:confusion_matrix}. This gave us an idea of what digits were easy or difficult to distinguish between. We then replicate Experiment 1 with a different couple of digits, namely $1$ and $8$, which are easier to distinguish than the digits $4$ and $9$. Results are depicted in Figure \ref{fig:ablation18}. We can draw the same conclusions of Experiment 1.
	
	%%%%%% BIBLIOGRAPHY, commands below break long URLs	
	\newcounter{biburlucpenalty}
	\newcounter{biburllcpenalty}
	\setcounter{biburllcpenalty}{7000}
	\setcounter{biburlucpenalty}{8000}
	\bibliography{references}
\end{document}